\documentclass{gtart}
\gtart

\usepackage{amsthm}
\usepackage{amsgen}
\usepackage{amsmath,amssymb}
\usepackage{curvesls}
\usepackage{epic}
\usepackage[all]{xy}
\usepackage[dvips]{graphicx}
\usepackage{psfrag}
\usepackage[T1]{fontenc}
\usepackage[sc]{mathpazo}
\newtheorem{thm}{Theorem}[section]
\newtheorem{lem}[thm]{Lemma}
\newtheorem{cor}[thm]{Corollary}
\newtheorem{prop}[thm]{Proposition}
\newtheorem{eg}[thm]{Example}
\newtheorem{prob}[thm]{Problem}

\theoremstyle{definition}
\newtheorem{defn}[thm]{Definition}

\newcommand{\Prime}{{\mathcal{P}}}
\newcommand{\Real}{{\mathbb R }}
\newcommand{\Rat}{{\mathbb Q}}

\newcommand{\Zed}{{\mathbb Z }}
\newcommand{\Nat}{{\mathbb N }}

\newcommand{\Diff}{{\mathrm{Diff}}}

\newcommand{\EK}[1]{{\mathrm{EC}({#1})}}
\newcommand{\ED}[1]{{\mathrm{ED}({#1})}}

\newcommand{\Cu}{{\mathcal{C}}}

\newcommand{\SD}{{\mathcal{SD}}} 
\newcommand{\SC}{{\mathcal{SC}}} 
\newcommand{\SP}{{\mathcal{SP}}} 

\newcommand{\K}{{\mathcal{K}}}

\newcommand{\upepsilon}{{ \xy
                   *{\xy (0,-3);(0,-2) **\dir{-} \endxy}
                   *\cir<3pt>{d^u}
                   \endxy}}

\begin{document}

\title{An operad for splicing}
\author{Ryan Budney}

\address{
Mathematics and Statistics, University of Victoria\\
PO BOX 3060 STN CSC, Victoria, B.C., Canada V8W 3R4
}
\email{rybu@uvic.ca}

\begin{abstract} 
A new topological operad is introduced, called the splicing operad.  This operad acts on a broad class of spaces of self-embeddings $N \to N$ where $N$ is a manifold.  The action of this operad on $\EK{j,M}$ (self embeddings $\Real^j \times M \to \Real^j \times M$ with support in $I^j\times M$) is an extension of the action of the operad of $(j+1)$-cubes on this space defined in \cite{Bud}. Moreover the action of the splicing operad encodes a version of Larry Siebenmann's \cite{BS, Sieb} splicing construction for knots in $S^3$ in the $j=1$, $M=D^2$ case, for which we denote the splicing operad $\SP_{3,1}$.  The space of long knots in $\Real^3$ (denoted $\K_{3,1}$) was shown to be a free algebra over the $2$-cubes operad with free generating subspace $\Prime \subset \K_{3,1}$, the subspace of long knots that are prime with respect to the connect-sum operation \cite{Bud}.  One of the main results of this paper is that $\K_{3,1}$ is free with respect to the splicing operad $\SP_{3,1}$ action, but the free generating space is the significantly smaller space of torus and hyperbolic knots $\mathcal{TH} \subset \K_{3,1}$.  Moreover, we show that $\SP_{3,1}$ is a free product of two operads.  The first free summand of $\SP_{3,1}$ is a semi-direct product $\Cu_2 \rtimes O_2$ operad which is {\it not} equivalent to the framed discs operad.  The second free summand of $\SP_{3,1}$ is a free $\Sigma \wr O_2$-operad, free on $\Sigma \wr O_2$-spaces which encode cabling and hyperbolic satellite operations, moreover the $\Sigma \wr O_2$-homotopy-type of these spaces is determined by finding adapted maximal symmetry positions for hyperbolic links in $S^3$.  This is an in-principle explicit description of the homotopy-type of the space of knots in $S^3$, and modulo the rather difficult problem of determining the symmetry groups of a class of hyperbolic links and their actions on the cusps, this is a closed form description of the homotopy-type. 
\end{abstract}

\maketitle

\section{Introduction}\label{intro}

In 1949 Schubert \cite{Sch1} proved that long knots in $\Real^3$ have a unique decomposition 
into prime knots.  A concrete statement of his theorem is that there is a homotopy-associative pairing 
$\K_{3,1} \times \K_{3,1} \to \K_{3,1}$ called the connect-sum operation which turns $\pi_0 \K_{3,1}$ (the isotopy 
classes of long knots) into a free commutative monoid.  The generators are called prime knots. 
The idea for why $\pi_0 \K_{3,1}$ is commutative is summarized in the diagram below. 

{
\psfrag{fg}[tl][tl][1][0]{$f \# g$}
\psfrag{gf}[tl][tl][1][0]{$g \# f$}
\psfrag{f}[tl][tl][1][0]{$f$}
\psfrag{g}[tl][tl][1][0]{$g$}
$$\includegraphics[width=12cm]{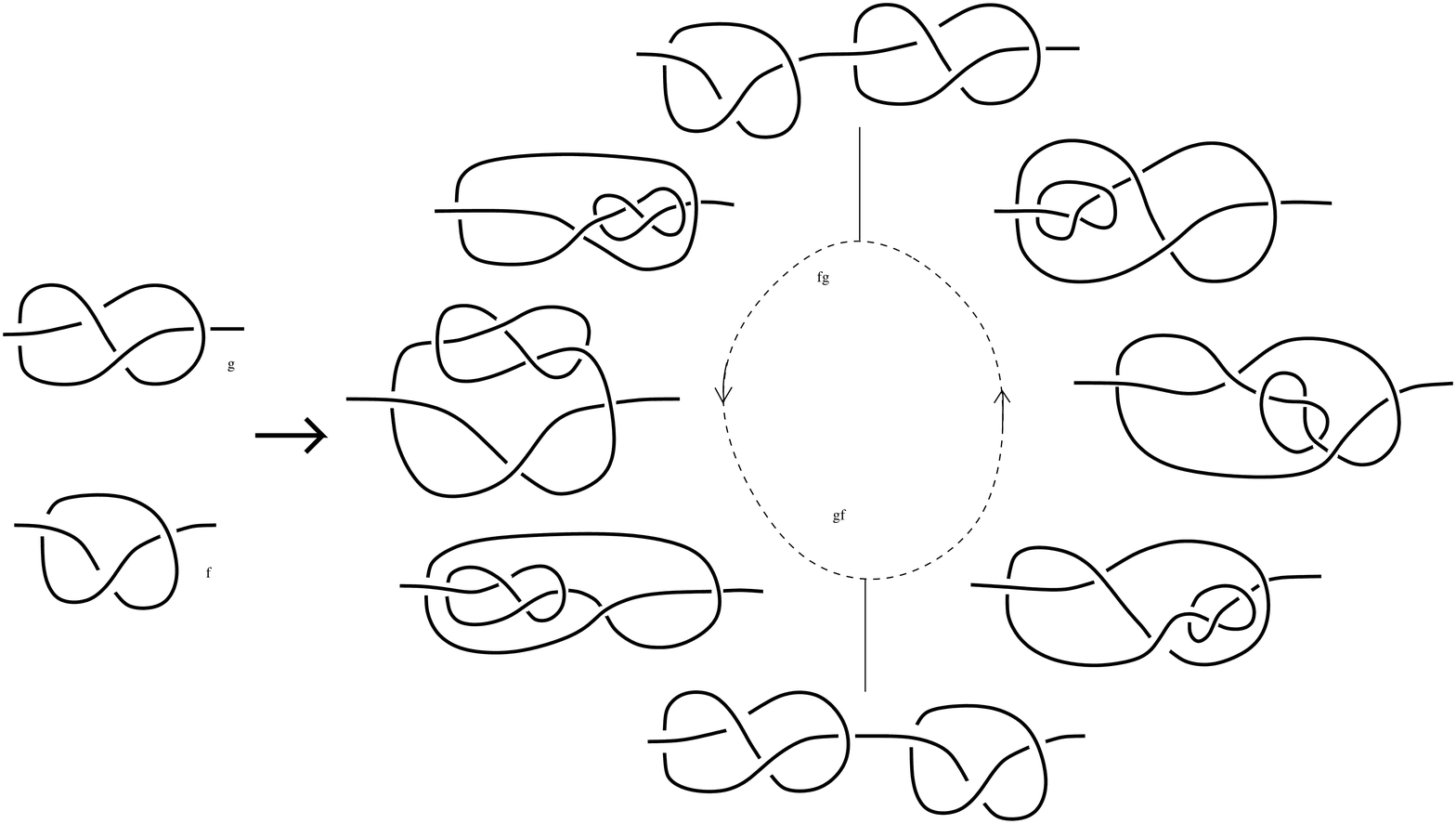}$$
}

`Little cubes and long knots' \cite{Bud} can be viewed as a space-level generalization of Schubert's work.  Schubert's homotopy-associative connect-sum mapping $\K_{3,1}\times \K_{3,1} \to \K_{3,1}$ is enhanced to an action of the operad of $2$-cubes $\Cu_2$ on $\K_{3,1}$, giving an explicit operadic parametrization of the kinds of isotopies depicted above.  The main theorem of \cite{Bud} is that $\K_{3,1}$ is free as an algebra over the $2$-cubes operad $\K_{3,1} \simeq \Cu_2(\Prime \sqcup \{*\})$, which when we apply $\pi_0$ recovers Schubert's result, since $\Cu_2(\Prime \sqcup\{*\})\simeq \sqcup_{n=0}^\infty \left(\Cu_2(n)\times_{\Sigma_n} \Prime^n\right)$. 

Schubert went on to further decompose knots using what he called {\it satellite operations} in his massive paper {\it Knoten und Vollringe} \cite{Sch2}.  As Schubert noticed, there are many ways to construct the same knot via distinct satellite operations.  In hindsight we know this was partially an accident of notation, as Schubert's notion of satellite operation was too linearly presented to see the symmetries inherent in the process of constructing satellite knots. Further, satellite constructions produce knots with incompressible tori in their complements, so the uniqueness statement must be tied to the JSJ-decomposition of 3-manifolds. The uniqueness statement for the JSJ decomposition is quite delicate and in some sense its delicate nature was a key factor in it being difficult to find.  It has been pointed out several times since and in several different contexts \cite{BS, EN, BudJSJ, Ikeda} that when reinterpreted via Larry Siebenmann's less linearly-ordered notion of {\it splicing} \cite{Sieb} there is a unique decomposition theorem for satellite knots.   

The primary point of this paper is to do for splicing what `little cubes and long knots' \cite{Bud} did for the connect-sum operation.  An operadic space-level encoding of splicing is given in Proposition \ref{def2}. Theorem \ref{freenessthm} shows $\K_{3,1}$ to be a free algebra over the splicing operad $\SP_{3,1}$, with free generating subspace the torus and hyperbolic knots $\mathcal{TH}$, i.e. $\K_{3,1} \simeq \SP_{3,1}(\mathcal{TH})$.  This provides a pleasant linkage between the low-dimensional topologists' view of knots (that torus and hyperbolic knots are in some sense the most essential), with the algebraic topologist's language of operads.   Further, it forms a link between the usage of trees in the study of operads to depict iterated composites of the structure maps with trees in 3-manifold theory, used to depict the structure of the JSJ-decomposition of a knot or link complement in $S^3$.  This is closely related to the somewhat unsatisfactory recursive structure of the homology of the long knot space $\K_{3,1}$ viewed as an algebra over the operad of $2$-cubes \cite{Cohen}.  The main result of \cite{Bud} is that $\K_{3,1}$ as an algebra over the operad of little $2$-cubes is free, where the free generating subspace is the space $\Prime \subset \K_{3,1}$ of prime long knots, i.e. $\K_{3,1} \simeq \Cu_2(\Prime \sqcup \{*\})$. As was observed in \cite{Cohen} and \cite{BudJSJ}, the homology of $\Prime$ has a deeper structure coming from the splicing decomposition of knots, forcing $H_*(\K_{3,1},\Rat)$ to reappear in shifted degrees inside $H_*(\Prime,\Rat)$ in many ways.  The non-operadic nature of the description of $\K_{3,1}$ given in \cite{KS} is non-uniform and somewhat frustrating. These complications largely disappear when $\K_{3,1}$ is viewed through the lens of the splicing operad $\SP_{3,1}$.  Theorem \ref{operadstruc} shows the splicing operad $\SP_{3,1}$ to be a free product (in the category of $\Sigma \wr O_2$-operads) of $\Cu_2 \rtimes O_2$ (not the framed discs operad, but a different semi-direct product) and various free operads.  The other free summands of $\SP_{3,1}$ correspond to cabling operations and hyperbolic satellite operations.  Moreover, all these summands with the sole exception of $\Cu_2 \rtimes O_2$ are free operads, freely generated on certain $\Sigma_k^* \wr O_2$-spaces, whose equivariant homotopy-type is identified in Theorem \ref{operadstruc}. 

A secondary point of this paper is that these techniques extend beyond the realm of classical knots. There are splicing operads that act on a wide class of spaces of self-embeddings $N \to N$, for $N$ a compact manifold.  This includes the spaces $\EK{j,M}$ and $\ED{j,M}$ \cite{Budsurv} of self-embeddings $\Real^j \times M \to \Real^j \times M$ with support contained in $[-1,1]^j \times M$ and $D^j \times M$ respectively, but the definition of the splicing operad applies to more general self-embedding spaces, some are discussed briefly in Section \ref{futuredir}.  In particular, the splicing operad for the `cubically supported embedding spaces' $\EK{j,M}$ is generally richer than the action of the corresponding action of the $(j+1)$-cubes operad on $\EK{j,M}$.  The splicing operad differs significantly from the operad of cubes, in that the splicing operad is an infinite-dimensional Frech\'et manifold, i.e. it is `big' when compared to many traditional operads, which tend to be levelwise finite-dimensional.   Another large-scale difference is that while the operad of $(j+1)$-cubes acts on the space $\EK{j,M}$ for all compact manifolds $M$, there are distinct splicing operads for $\EK{j,M}$ and $\EK{j,N}$ provided $M$ and $N$ are distinct.  Perhaps this new operad will lead to new insights into the homotopy-types of these embedding spaces. 

This paper was influenced by conversations with Jim McClure, Paolo Salvatore and Allen Hatcher.  Thanks to BIRS hosting Allen Hatcher's 65th birthday party where I had the opportunity to run these ideas past the participants. Thanks to the University of Rome, Tor Vergata, for hosting me in the summer of 2009 where these ideas indirectly started fermenting. Thanks also to Toshitake Kohno, the University of Tokyo and the Institute for the Physics and the Mathematics of the Universe (IPMU) for hosting me in the winter of 2010 and 2012.   Thanks to Victor Turchin and Tom Goodwillie for comments on the initial drafts of this manuscript.

\section{The operad of overlapping $n$-cubes}

The point of this section is to provide a motivating result, vaguely this is a `flattening' of the operad of little $(n+1)$-cubes to an equivalent operad called the operad of overlapping $n$-cubes. None of the main results of this paper depend significantly on this section. These results are provided as context, as part of the train of thought leading up to the construction in Section \ref{splicing}, which might otherwise seem as uninspired. The point of this construction is that the operad of overlapping $n$-cubes has a more natural action on embeddings spaces, equivalent the the action of the operad of little $(n+1)$-cubes on $\EK{n,M}$. 

\begin{defn}\label{operadDEF} 
A topological $\Sigma$-operad is a collection of right $\Sigma_n$-spaces $\mathcal O(n)$ for $n \in \{0,1,2,\cdots\}$ and maps $$\mathcal O(k) \times \left( \mathcal O(j_1) \times \cdots \times \mathcal O(j_k) \right) \to \mathcal O(j_1 + \cdots + j_k)$$ satisfying an (1) associativity, (2) symmetry and (3) identity axiom.  Given $J \in \mathcal O(k)$ and $L_i \in \mathcal O(j_i)$ for $i \in \{1,2,\cdots,k\}$ denote the image of $(J,L_1,\cdots,L_k)$ under the above map by $J.L$.  
\begin{enumerate}
\item The associativity condition is that $J.(L.M) = (J.L).M$ whenever this makes sense, i.e. $M = (M_{1,1}, \cdots, M_{1,j_1}, M_{2,1}, \cdots, M_{2,j_2}, \cdots, M_{k,1}, \cdots, M_{k,j_k})$ with each $M_{a,j_b}$ belonging to the operad $\mathcal O=\sqcup_{n=0}^\infty \mathcal O(n)$. 
\item The symmetry axiom is that $(J.\sigma).(\sigma^{-1}.L) = (J.L).\overline{\sigma}$.  
We interpret $L$ as a $k$-tuple $L=(L_1,\cdots,L_k)$, so the left action of $\sigma^{-1}$ on $L$ is
$\sigma^{-1}.L = (L_{\sigma(1)}, \cdots, L_{\sigma(k)})$.  $\overline{\sigma} \in \Sigma_{j_1+\cdots+j_k}$
is the associated block permutation to $\sigma$.  Similarly there is a symmetry condition
$(J.L).\theta = J.(L.\theta)$ provided $\theta = \theta_1 \times \cdots \times \theta_k$ with $\theta_i \in \Sigma_{j_i}$
for all $i\in \{1,2,\cdots, k\}$. 
\item The identity axiom is that there is an element $I \in \mathcal O(1)$ such that $I.L = L$ for all $L \in \mathcal O$, and that $J.(I, \cdots, I) = J$ for all $J \in \mathcal O$. 
\end{enumerate} 
An action of the operad $\mathcal O$ on a space $X$ is a sequence of maps $\mathcal O(n) \times X^n \to X$
for $n \in \{0,1,\cdots\}$ satisfying an (1) associativity, (2) symmetry and (3) identity axiom. As above, let $J \in \mathcal O(k)$, and $L_i \in \mathcal O(j_i)$ for $i\in \{1,2,\cdots,k\}$.  
\begin{enumerate}
\item The associativity condition demands that $(J.L).x = J.(L.x)$ provided $x \in X^{j_1 + \cdots + j_k}$
\item The symmetry condition demands that $(J.\sigma).x = J.(\sigma.x)$ where the left action of $\sigma$ on 
$X^k$ is given by $\sigma.(x_1,\cdots,x_k) = (x_{\sigma^{-1}(1)},\cdots,x_{\sigma^{-1}(k)})$. 
\item The identity condition demands that if $I \in \mathcal O(1)$ is the identity of $\mathcal O$, then $I.x = x$ for all $x \in X$. 
\end{enumerate}
\end{defn}

Operads were originally designed as a category theoretic analogue of universal algebras.  The above definition immediately generalizes to operads in symmetric monoidal categories, see \cite{MSS, May} for example.   The space $\mathcal O(0)$ will be called the {\it base} of the operad (sometimes called the $0$-th operadic grading, or the constants of the operad). Notice that the structure maps of $\mathcal O$ restrict to an action of $\mathcal O$ on the base $\mathcal O(k) \times \left( \mathcal O(0) \times \cdots \times \mathcal O(0)\right) \to \mathcal O(0)$, this will be called the {\it augmentation action}.  If the base consists of a single point, the operad is said to be {\it unitial}. If $\mathcal O(1)$ consists of a single point the operad is said to be {\it reduced}.  Some authors include as part of their definition that the base is empty \cite{MSS, May}, although this is not a uniform requirement among authors.  The first operad discussed in the literature is the operad of little cubes, which appears with both an unbased and unitial variant.  In this paper the cubes operad (Definition \ref{overlappingcubes}) is unitial.  For an operad with non-empty base the structure maps give {\it degeneracy maps} $\mathcal O(n) \times \left( \mathcal O(1)^i \times \mathcal O(0) \times \mathcal O(1)^{n-i-1} \right) \to \mathcal O(n-1)$, which when restricted to $(I, I, \cdots, I, *, I, \cdots, I) \in \mathcal O(1)^i \times \mathcal O(0) \times \mathcal O(1)^{n-i-1}$ gives maps $\mathcal O(n) \to \mathcal O(n-1)$, here $* \in \mathcal O(0)$ is a choice of base-point. 

\begin{defn}\label{overlappingcubes}
An increasing affine-linear function $[-1,1] \to [-1,1]$ is a {\it little interval}.  A product of little intervals
$[-1,1]^n \to [-1,1]^n$ is a {\it little $n$-cube}.  The space $\Cu_n(j)$ is the collection of $j$-tuples of
little $n$-cubes whose images are required to have disjoint interiors, $\Cu_n(0)=\{*\}$ is the empty cube.  The collection 
$\Cu_n = \sqcup_{j=0}^\infty \Cu_n(j)$ is the operad of little $n$-cubes, it is a $\Sigma$-operad with 
structure maps $\Cu_n(k) \times \left( \Cu_n(j_1) \times \cdots \times \Cu_n(j_k) \right) \to \Cu_n(j_1+\cdots+j_k)$
defined by $(L,J_1,\cdots,J_k) \longmapsto (L_1 \circ J_1, \cdots, L_k \circ J_k)$ 
and $\Cu_n(j) \times \Sigma_j \to \Cu_n(j)$ given by $(L, \sigma) \longmapsto L\circ \sigma$.  We take $\Sigma_j = \mathrm{Aut} \{ 1,2, \cdots, j \}$ throughout the paper. Sometimes we will further think of $\Sigma_j$ as the subgroup of $\mathrm{Aut} \{0,1,2,\cdots,j\}$ that fix $0$, but in this case $\Sigma_j$ will be denoted $\Sigma_j^*$. 

A collection of {\it $j$ overlapping $n$-cubes} is an equivalence class of pairs $(L, \sigma)$ where  $L=(L_1,\cdots,L_j)$, each $L_i$ is a little $n$-cube and $\sigma \in \Sigma_j$.  Two collections of $j$ overlapping $n$-cubes  $(L,\sigma)$ and $(L',\sigma')$ are taken to be {\it equivalent} provided $L = L'$ and whenever the interiors of $L_i$ and $L_k$ intersect $\sigma^{-1}(i) < \sigma^{-1}(k) \Longleftrightarrow 
\sigma'^{-1}(i) < \sigma'^{-1}(k)$.  Given $j$ overlapping $n$-cubes $(L_1,\cdots,L_j,\sigma)$ we say the $i$-th cube $L_i$ is at {\it height} $\sigma^{-1}(i)$. $\sigma(1)$ is the index of the {\it bottom} cube, and $\sigma(j)$ is the index of the {\it top} cube.  Let $\Cu_n'(j)$ be the space of all $j$ overlapping $n$-cubes, with the quotient topology induced by the equivalence relation. 

The structure map 
$$\Cu_n'(k) \times \left( \Cu_n'(j_1) \times \cdots \times \Cu_n'(j_k) \right) \to \Cu_n'(j_1 + \cdots + j_k)$$
is defined by
$$\left((L,\sigma), (J_1,\alpha_1), \cdots, (J_k, \alpha_k)\right) \longmapsto 
  ((L_1\circ J_1, \cdots, L_k\circ J_k), \beta)$$
the permutation $\beta$ is given for $1 \leq a \leq k$, $1 \leq b \leq j_a$
$$\beta^{-1}\left(\sum_{i<a} j_i + b\right) = \left( \sum_{i<\sigma^{-1}(a)} j_{\sigma(i)}\right) +\alpha^{-1}_a(b).$$
This permutation is obtained by taking the lexicographical order on the set
$\{(a,b) : a \in \{1,\cdots,k\}, b \in \{1,\cdots,j_a\}\}$ and then identifying with 
$\{1, 2, \cdots,j_1+\cdots+j_k\}$ in the order-preserving way. 
\end{defn}

Next we will adapt the action of $\Cu_{j+1}$ on $\EK{j,M}$ to be an action of $\Cu_j'$ on $\EK{j,M}$. 
First a reminder of the definition and geometric context for the action of $\Cu_{j+1}$ on $\EK{j,M}$. 

\begin{defn}\label{longdef} A (thin) long knot is a smooth embedding
$\Real^j \to \Real^n$ which agrees with the standard embedding
$x \longmapsto (x,0)$ outside of the cube $I^j = [-1,1]^j$. The space of
thin long knots is denoted $\K_{n,j}$.  In various situations one might want
to replace $I^j$ in this definition by $D^j = \{ x \in \Real^j : |x| \leq 1 \}$. 
We distinguish between these definitions by saying the knot has {\it cubical} 
support versus being {\it supported on a disc}.  It's an elementary rescaling
argument that the inclusion $\K_{n,j}^{\text{disc}} \to \K_{n,j}^{\text{cubical}}$
is a homotopy-equivalence. 

{
\psfrag{ek1d2}[tl][tl][1.2][0]{$f \in \EK{1,D^2}$}
\psfrag{p1}[tl][tl][0.8][0]{$1$}
\psfrag{m1}[tl][tl][0.8][0]{$-1$}
$$\includegraphics[width=10cm]{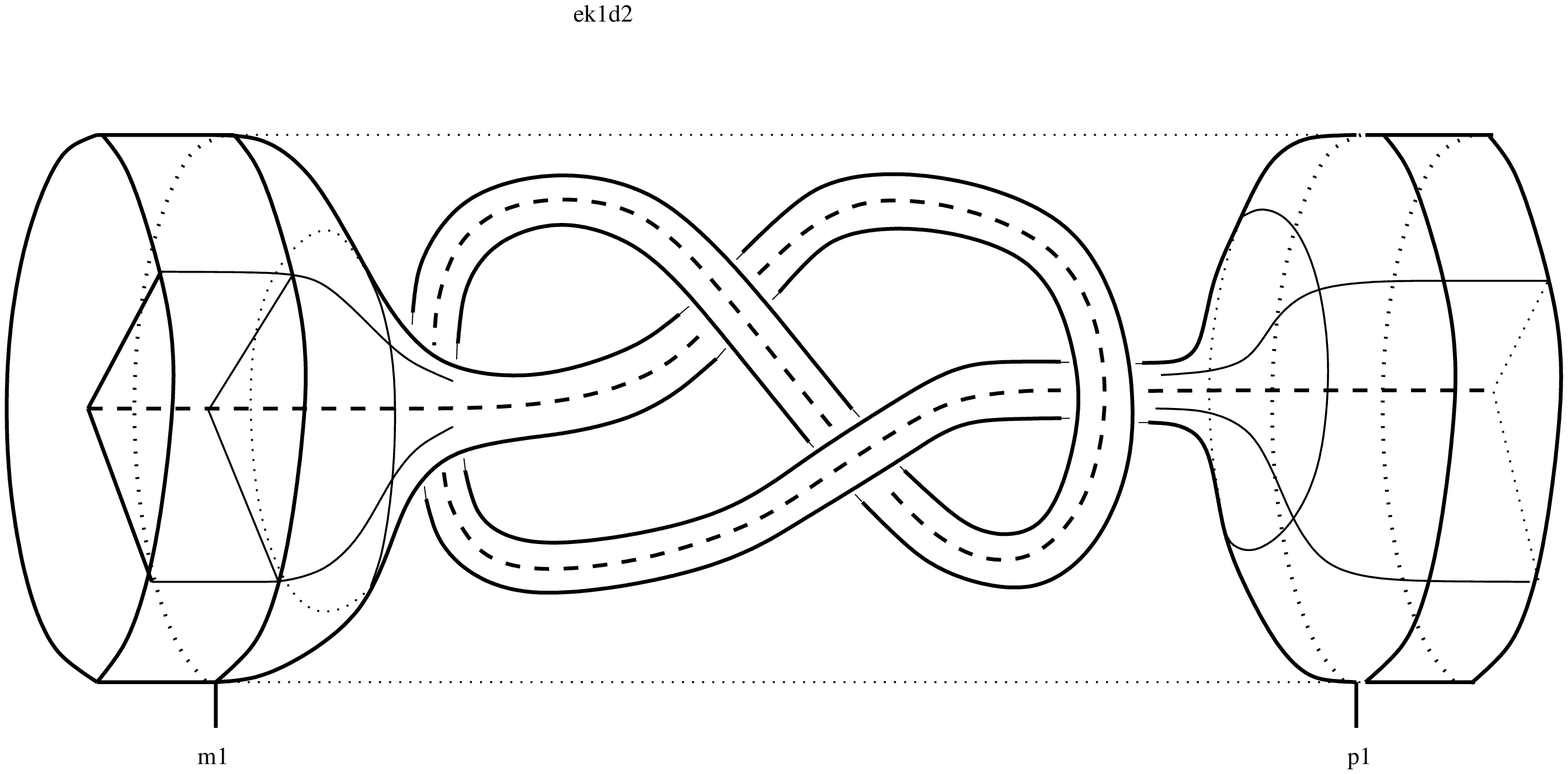}$$
}

A (fat) long knot is an embedding $f : \Real^j \times M \to \Real^j \times M$
such that $supp(f) \subset I^j \times M$.  The space of fat long knots is
denoted $\EK{j,M}$.  The restriction map $\EK{j,D^{n-j}} \to \K_{n,j}$ given by
$f \longmapsto f_{|\Real^j \times \{0\}}$ is a fibration whose fibre has the
homotopy-type of $\Omega^j SO_{n-j}$. So typically $\EK{j,D^{n-j}}$ is called the space of {\it framed}
long knots, as it consists of knots together with an explicit trivialization of a tubular
neighbourhood. The notation $\mathrm{EC}$ is meant to indicate `embeddings with
cubical support.'  $\EK{1,D^2}$ has the homotopy-type of $\K_{3,1} \times \Zed$ since
the fibration $\EK{1,D^2} \to \K_{3,1}$ splits at the fibre, with splitting given by
the linking-number of $f_{|\Real \times \{(0,0)\}}$ and $f_{|\Real \times \{(1,0)\}}$. 
Thus $\K_{3,1}$ has the homotopy-type of $\hat \K_{3,1} \subset \EK{1,D^2}$
and $\EK{1,D^2} = \Zed \times \hat \K_{3,1}$, where $\hat \K_{3,1}$ is the subspace of
$\EK{1,D^2}$ consisting of knots $f$ where the above linking number is zero.  The homotopy-equivalence 
$\hat \K_{3,1} \to \K_{3,1}$ is the restriction map \cite{Bud}.   As with long knots, 
if one replaces every occurrence of $I^j$ by $D^j$ one gets a homotopy-equivalent space 
$\ED{j,M}$, the inclusion $\ED{j,M} \to \EK{j,M}$ being a homotopy-equivalence.
\end{defn}

The choice of usage of discs or cubes in the definitions of $\K_{n,j}$, $\ED{j,M}$ and $\EK{j,M}$
becomes important when one wants to study group actions on these spaces.  For example, 
$\K_{n,j}^{\text{disc}}$ admits an action of $O_j$ (by conjugation), while 
$\K_{n,j}^{\text{cubical}}$ does not.  Further, the family of spaces $\K_{n,j}^{\text{cubical}}$ fits into 
a pseudoisotopy fibration sequence (see \cite{Budsurv}), while the family $\K_{n,j}^{\text{disc}}$ does not. 

We assemble the ingredients of the action of $\Cu_{j+1}$ on $\EK{j,M}$.  Given a little
$j$-cube $L$ and $f \in \EK{j,M}$ the {\it rescaling} of $f$ by $L$ is  $L.f = (L \times Id_M) \circ f \circ (L \times Id_M)^{-1}$. For this to make sense, reinterpret $L$ as its unique affine-linear extension $L : \Real^j \to \Real^j$. 
Given a $(j+1)$-cube $L$, write it as a product $L^\pi \times L^\nu$ where $L^\pi$ is a $j$-cube and 
$L^\nu$ is a $1$-cube. Let $L^t = L^{\nu}(-1)$.  Given $n$ little $(j+1)$-cubes, 
$L=(L_1,\cdots, L_n)\in \Cu_{j+1}(n)$ define the $n$-tuple of (non-disjoint) little $j$-cubes 
$L^\pi = (L_1^\pi,\cdots, L_n^\pi)$. Similarly define $L^t \in I^j$ 
by $L^t=(L_1^t,\cdots, L_n^t)$. The action of $\Cu_{j+1}$ on $\EK{j,M}$ \cite{Bud} was defined as  
$\kappa_n : \Cu_{j+1}(n) \times \EK{j,M}^n \to \EK{j,M}$ for $n \in \{1,2,\cdots\}$ which is given by
$$\kappa_n(L_1,\cdots,L_n,f_1,\cdots,f_n) =
  L^\pi_{\sigma(n)}.f_{\sigma(n)}\circ L^\pi_{\sigma(n-1)}.f_{\sigma(n-1)}\circ\cdots \circ L^\pi_{\sigma(1)}.f_{\sigma(1)}$$
where $\sigma : \{1,\cdots, n\} \to \{1,\cdots,n\}$ is any permutation
such that $L^t_{\sigma(n)} \geq L^t_{\sigma(n-1)} \geq \cdots \geq L^t_{\sigma(1)}$.  
Notice that the action of $\Cu_{j+1}$ on $\EK{j,M}$ has a rather coarse dependence on the cubes $L$, in 
that only the relative ordering specified by $\sigma$ matters, much of the information given by $L^\nu$ 
is irrelevant.  This will be made precise in Proposition \ref{flataction}. 

{
\psfrag{p1}[tl][tl][0.6][0]{$-1$}
\psfrag{m1}[tl][tl][0.6][0]{$1$}
\psfrag{mu}[tl][tl][1][0]{}
\psfrag{e}[tl][tl][1][0]{$f$}
\psfrag{f}[tl][tl][0.8][0]{$L$}
\psfrag{e.f}[tl][tl][0.9][0]{$L.f$}
\psfrag{,}[tl][tl][1][0]{$,$}
\psfrag{t}[tl][tl][1][0]{$\{0\}^j\times \Real$}
\psfrag{pit}[tl][tl][0.9][0]{$L^t$}
\psfrag{pic}[tl][tl][0.9][0]{$L^\pi$}
\psfrag{rn}[tl][tl][1][0]{$\Real^j\times\{0\}$}
$$\includegraphics[width=8cm]{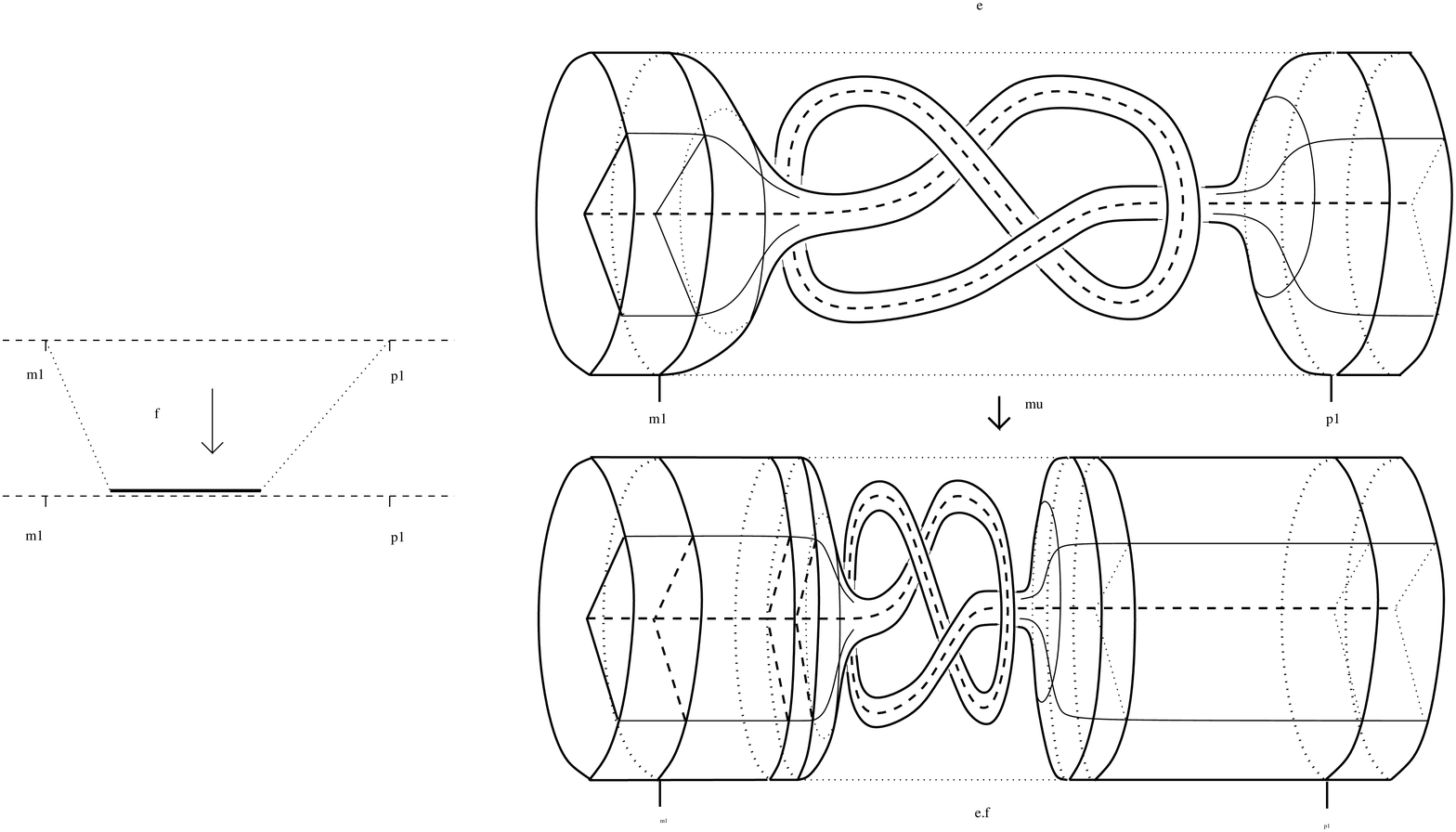} \hskip 1cm \includegraphics[width=5cm]{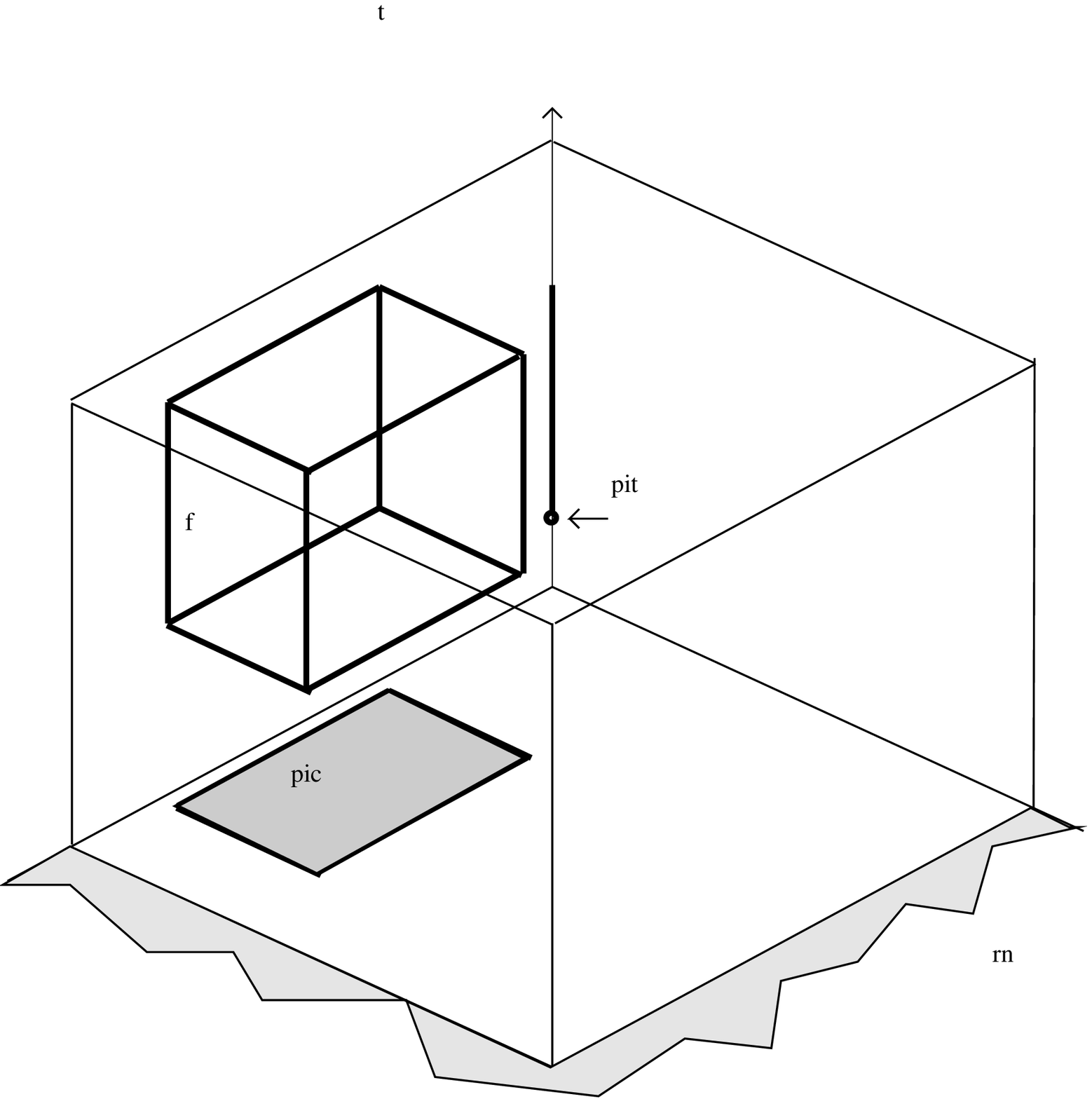}$$
\centerline{Rescaling $f$ by $L$ in $\EK{1,D^2}$ \hskip 3cm Projection $L \longmapsto L^\pi$ and $L^t$}
}

\begin{eg}
{
\psfrag{1}[tl][tl][1][0]{$L_1$}
\psfrag{2}[tl][tl][1][0]{$L_2$}
\psfrag{p1}[tl][tl][0.7][0]{$1$}
\psfrag{m1}[tl][tl][0.7][0]{$-1$}
\psfrag{l1t}[tl][tl][1][0]{$L_1^t$}
\psfrag{l2t}[tl][tl][1][0]{$L_2^t$}
\psfrag{comma}[tl][tl][1][0]{,}
\psfrag{nu3}[tl][tl][1][0]{$\kappa_2$}
\psfrag{f1}[tl][tl][1][0]{$f_1$}
\psfrag{f2}[tl][tl][1][0]{$f_2$}
$$\includegraphics[width=15cm]{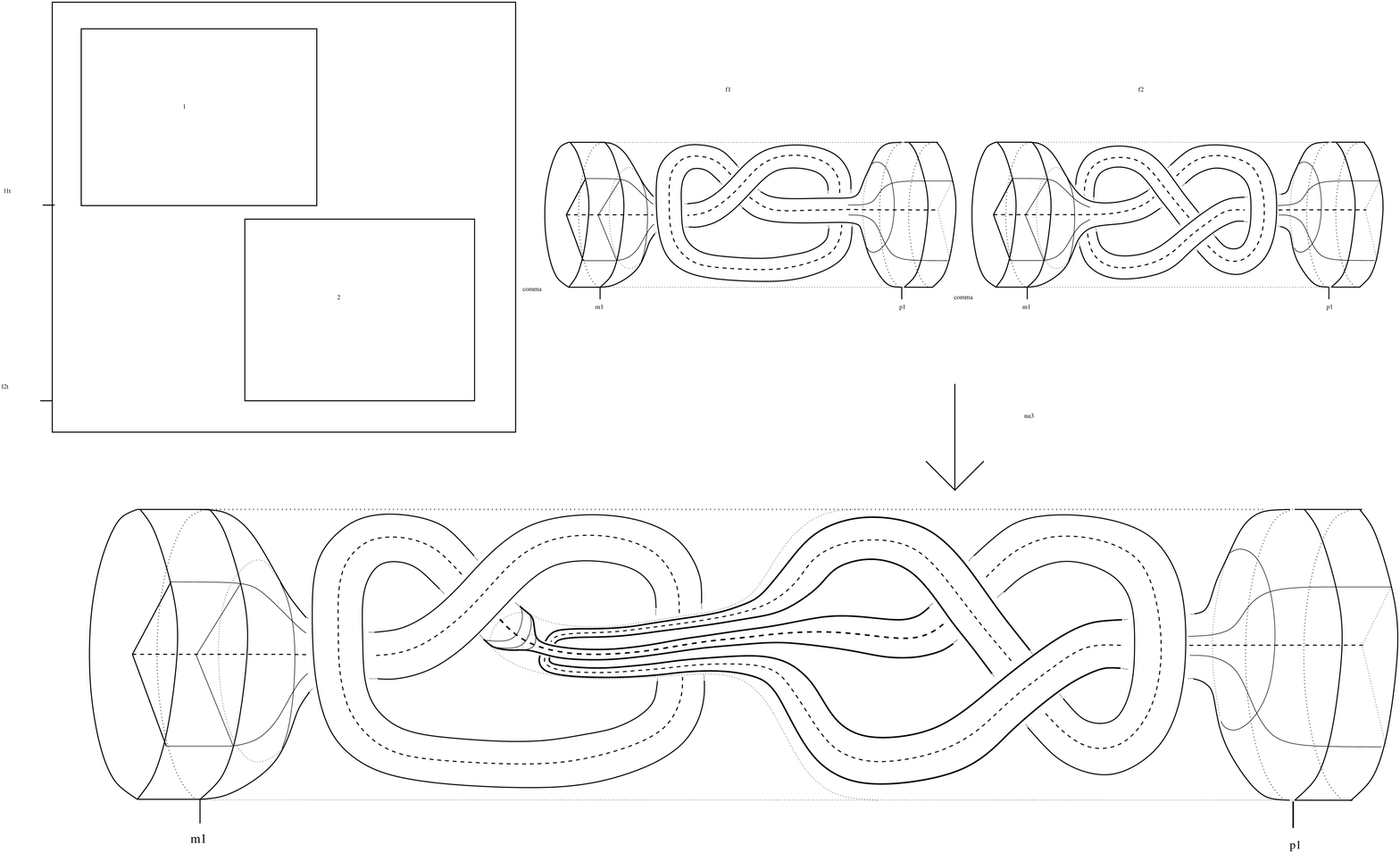}$$
\centerline{$L^t_1 > L^t_2$ so $\sigma = (12)$ and $\kappa_2(L_1,L_2,f_1,f_2)=L^\pi_1.f_1 \circ L^\pi_2.f_2$}
}
\end{eg}

\begin{eg}
{
\psfrag{1}[tl][tl][1][0]{$L_1$}
\psfrag{2}[tl][tl][0.7][0]{$L_2$}
\psfrag{3}[tl][tl][1][0]{$L_3$}
\psfrag{p1}[tl][tl][0.7][0]{$1$}
\psfrag{m1}[tl][tl][0.7][0]{$-1$}
\psfrag{comma}[tl][tl][1][0]{,}
\psfrag{nu3}[tl][tl][1][0]{$\kappa_3$}
\psfrag{l3t}[tl][tl][1][0]{$L_3^t$}
\psfrag{l2t}[tl][tl][1][0]{$L_1^t$}
\psfrag{l1t}[tl][tl][1][0]{$L_2^t$}
\psfrag{f3}[tl][tl][1][0]{$f_3$}
\psfrag{phi1}[tl][tl][1][0]{$f_1$}
\psfrag{phi2}[tl][tl][1][0]{$f_2$}
\psfrag{phi3}[tl][tl][1][0]{$f_3$}
$$\includegraphics[width=15cm]{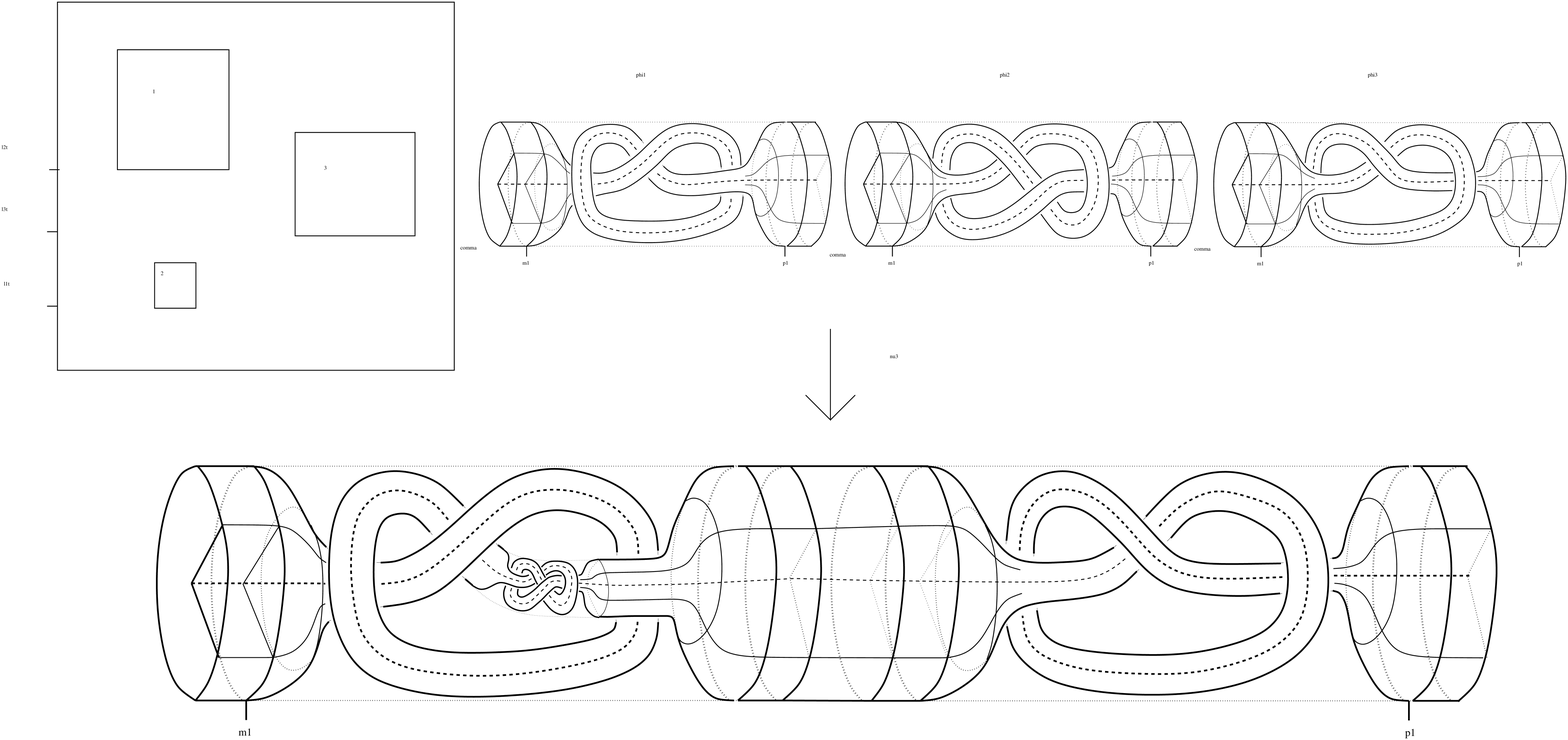}$$
\centerline{$L^t_1 > L^t_3 > L^t_2$ so $\sigma=(23)$ and $\kappa_3(L_1,L_2,L_3,f_1,f_2,f_3)=L^\pi_1.f_1\circ L^\pi_3.f_3 \circ L^\pi_2.f_2$}
}
\end{eg}

\begin{prop}\label{flataction}
$\Cu_j'$ is a multiplicative $\Sigma$-operad and the projection map $\Cu_{j+1} \to \Cu_j'$ given by $(L_1,\cdots,L_n) \longmapsto (L_1^\pi,\cdots,L_n^\pi, \sigma)$ as defined above is an operad map which is also a homotopy equivalence.
The maps $\kappa_n' : \Cu_j'(n) \times \EK{j,M}^n \to \EK{j,M}$ given by
$$\kappa_n'((L_1,\cdots,L_n,\sigma),(f_1,\cdots,f_n)) = L_{\sigma(n)}.f_{\sigma(n)} \circ \cdots \circ L_{\sigma(1)}.f_{\sigma(1)}$$ 
define an action of the operad $\Cu_j'$ on $\EK{j,M}$, and there is a commutative diagram
$$ \xymatrix{ \Cu_{j+1}(n) \times \EK{j,M}^n \ar[r]^-{\kappa_n} \ar[d] & \EK{j,M} \\ 
              \Cu_j'(n) \times \EK{j,M}^n \ar[ur]_-{\kappa'_n} & .} $$
\begin{proof}
To show $\Cu_j'$ is an operad, that $\kappa'$ is an action of the 
operad on $\EK{j,M}$ and that the above diagram commutes is mechanical, compare to the proof
of Theorem 5 in \cite{Bud}.  To see that the projection map $\Cu_{j+1}(n) \to \Cu_j'(n)$ is 
a homotopy-equivalence, notice that the fibre over any point in $\Cu_j'(n)$ is a convex polyhedron, 
the affine structure being given by the top and bottom coordinates of $L^\nu$.  The statement
that $\Cu_j'$ is a multiplicative operad means that $\Cu_j'$ contains the associative operad as
a sub-operad. This is elementary, as
$\{ (Id_{[-1,1]^j}, \cdots, Id_{[-1,1]^j}, Id_{\{1,\cdots,k\}}) : k\in \Nat \} \subset \Cu_j'$ is isomorphic to the associative operad. 
\end{proof}
\end{prop}

There are `overlapping' variants of operads of balls, operads of framed discs and the
operads of conformal balls \cite{fcube}.  For example, the operad of overlapping $n$-balls 
is equivalent to the operad of $(n+1)$-balls, but is also multiplicative.  The operad of
overlapping conformal $n$-balls is cyclic and multiplicative but it is not equivalent to the operad of
conformal $(n+1)$-balls.  It fibers over the operad of overlapping $n$-balls but the fibre 
consists of products of $SO_n$. 

\section{Operadic splicing}\label{splicing}

For knots in $S^3$, splicing has a particularly physical nature.  Splicing's role is to create new knots from old.  If a knot is sitting in front of you, with your hands reach out and `grab' the knot.  In this grabbed position, each hand forms a loop around a collection of strands of the knot.  In abstract, we represent this `grabbed position' by a knot together with a disjoint trivial link (it would be a $2$-component trivial link in the case of a single $2$-handed person grabbing the knot).  The second step involves isolating the strands grasped inside an individual hand, and performing a local modification on the knot.  The rough idea for how to perform the local modification is to cut the strands that pass through an individual  hand, and perform a local knotting operation on those loose ends, before re-gluing the strands together.  The important aspect of this heuristic is that splicing involves two steps, (1) the `grabbing' of the knot, represented in Definition \ref{splicedef} by a {\it knot generating link} (KGL) and (2) the local operation on the `grabbed' knot, which is Definition \ref{def1}, the splicing operation.

The notion of `splicing' was first described by Siebenmann \cite{Sieb} in his work on the JSJ-decomp\-ositions of homology spheres.  Splicing has its roots in Schubert's satellite operations \cite{Sch2}, but only came to prominence with the JSJ-decomposition of 3-manifolds.  In 1987 Bonahon and Siebenmann went on to explain splicing for knots and links in 3-manifolds in some detail, together with the JSJ-decomposition of the $\Zed_2$-cyclic branched cover of links in 3-manifolds \cite{BS} although their preprint has been out of distribution until recently. Eisenbud and Neumann's book \cite{EN} describes the splice decomposition of graph homology spheres in detail.   The refinement of splicing adapted specifically to knots and links in $S^3$ was given in \cite{BudJSJ}, of which some elements are sketched in this section. The main point of this section is the construction of an operad $\SC_j^{M}$ which acts on $\EK{j,M}$ (and $\SD_j^M$ acting on $\ED{j,M}$ respectively) for which the $M=D^2$ and $j=1$ case the operad's action is splicing in the sense of \cite{BudJSJ}, while it is closely related to splicing in the senses of \cite{BS, Sieb, EN}.  Section \ref{futuredir} sketches some further generalizations of these operads.

\begin{defn}\label{splicedef}
A knot-generating link (KGL) \cite{Bud} is an $(n+1)$-tuple $(L_0, L_1, \cdots, L_n)$ where $L_0 \in \K_{3,1}$ is a thin long knot, $L_i : S^1 \to [-1,1] \times D^2$ is an embedding for $i \in\{1,2,\cdots,n\}$ such that $(L_0, L_1, \cdots, L_n)$ are disjoint and $\{L_1, \cdots, L_n\}$ represents the $n$-component unlink.  We require $n$ to be non-negative $n \in \{0,1,2,3,\cdots\}$.
 
A {\it splicing diagram} is an enhanced or `fattened' KGL, allowing for a canonical definition of splicing. While KGL's were developed for the embedding space $\K_{3,1}$ \cite{Bud}, splicing diagrams will make sense for any embedding space of the form $\EK{j,M}$ or $\ED{j,M}$.  A splicing diagram for $\EK{j,M}$ is an equivalence class of $(n+2)$-tuple $(L_0, L_1, \cdots, L_n, \sigma)$ where $\sigma \in \Sigma_n$ is a permutation, $L_0 \in \EK{j,M}$, and $L_i : [-1,1]^j \times M \to [-1,1]^j \times M$ is an embedding for all $i \in \{1,2,\cdots,n\}$. The equivalence relation is given by $(L,\sigma) \sim (L',\sigma') \Longleftrightarrow L=L'$ together
with the relation that if $L_i( ([-1,1]^j)^\circ \times M) \cap L_j(([-1,1]^j)^\circ \times M) \neq \emptyset$ then $\sigma^{-1}(i) < \sigma^{-1}(j) \Longleftrightarrow \sigma'^{-1}(i)< \sigma'^{-1}(j)$, where $i,j \in \{1,2,\cdots,n\}$. There is a further
{\it continuity constraint} on a splicing diagram, that whenever $0 \leq \sigma^{-1}(i)<\sigma^{-1}(k)$, we require $\overline{L_{i}([-1,1]^j \times M)) \setminus L_{k}\left([-1,1]^j \times M\right)} \cap L_{k}\left(([-1,1]^j)^\circ \times \partial M\right) = \emptyset$, for any $i,j \in \{0,1,\cdots,n\}$. For the purposes of the continuity constraint, we use the convention $\Sigma_n \equiv \Sigma_n^* \subset \mathrm{Aut}\{0,1,\cdots,n\}$ (i.e. every $\sigma \in \Sigma^*_n$ 
satisfies $\sigma(0)=0$). Let $\SC_j^M(n) = \{ (L_0,L_1, \cdots, L_n, \sigma ) : \text{ is a splicing diagram} \}$, with the quotient topology induced by the equivalence relation $\sim$. Above we use the convention that if $X$ is a manifold with boundary $X^\circ$ denotes the interior $X^\circ = X \setminus \partial X$.
\end{defn}

Comments on choices made in the above definition: 
\begin{itemize}
\item[1)] If one wants to avoid manifolds-with-corners in the definition of a splicing diagram (as in Definition \ref{longdef}) replace all occurrences of $[-1,1]^j$ in Definition \ref{splicedef} with $D^j$, similarly replace $\EK{j,M}$ by $\ED{j,M}$.   There are situations in which either formalism appears to be the more appropriate, cubes for pseudo-isotopy fibrations \cite{Budsurv} and discs when interested in symmetry. Let $\SD_j^M$ denote the splicing operad using the discs formalism.  Notice this makes no difference in the $j=1$ case, i.e. $\SD_1^M = \SC_1^M$ and $\EK{1,M} = \ED{1,M}$ always. 

\begin{eg} A splicing diagram.
{
\psfrag{ek1d2}[tl][tl][1.1][0]{$(L,\sigma) \in \SD_1^{D^2}(3)$}
\psfrag{p1}[tl][tl][0.8][0]{$1$}
\psfrag{m1}[tl][tl][0.8][0]{$-1$}
\psfrag{l0}[tl][tl][1.1][0]{$L_0$}
\psfrag{l1}[tl][tl][0.8][0]{$L_1$}
\psfrag{l2}[tl][tl][0.8][0]{$L_2$}
\psfrag{l3}[tl][tl][0.8][0]{$L_3$}
$$\includegraphics[width=10cm]{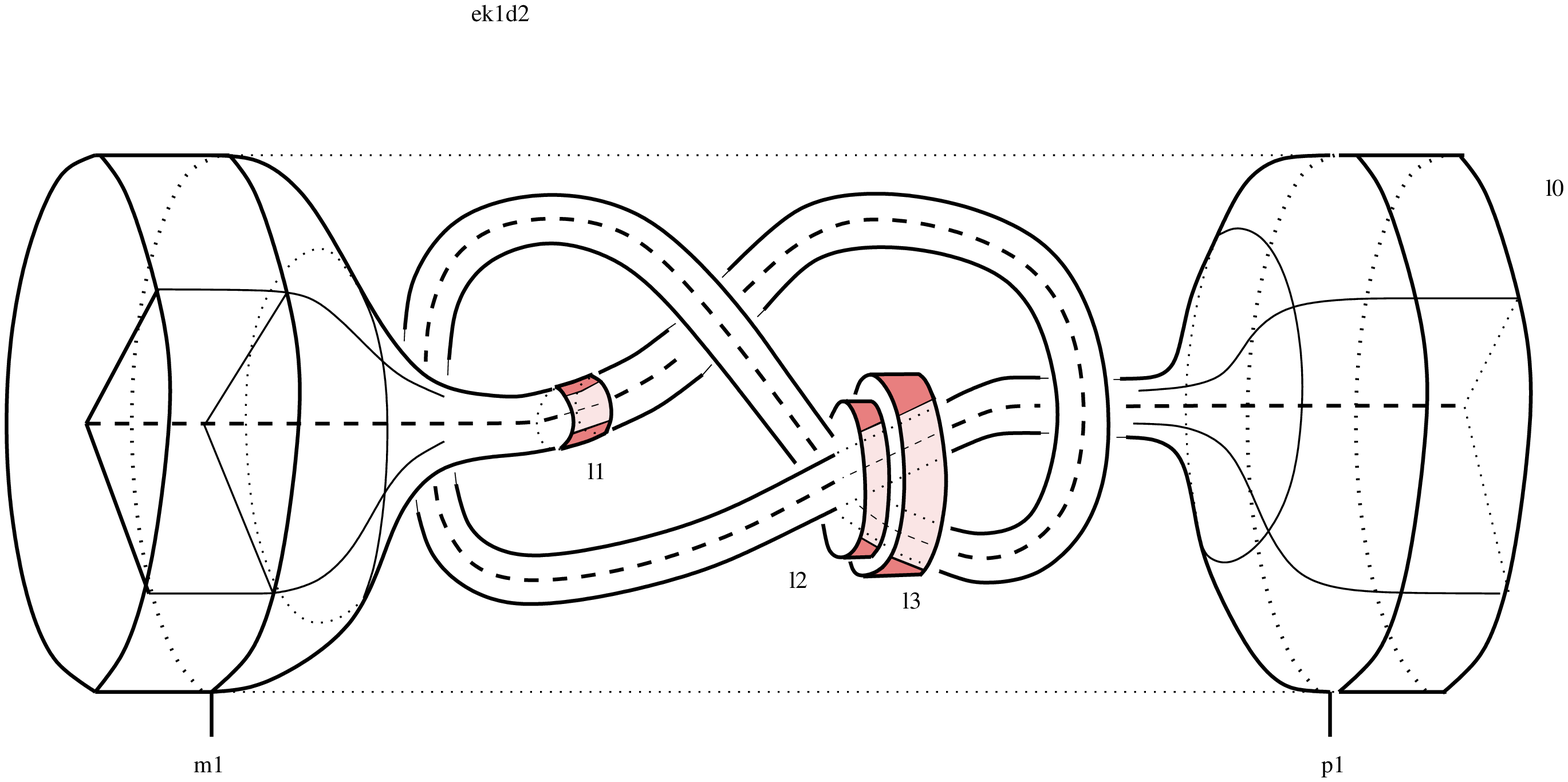}$$
Notice that $\sigma^{-1}(2) < \sigma^{-1}(3)$ is the only restriction on $\sigma \in \Sigma_3$
coming from Definition \ref{splicedef}, since the image of $L_2$ is partially contained in the image
of $L_3$. The order can't be reversed since $L_2([-1,1]\times S^1)$ intersects the interior of the
image of $L_3$. 
}
\end{eg}

\item[2)] To make sense of the continuity constraint some terminology is useful.  
Given an element $L = (L_0,L_1,\cdots,L_n,\sigma) \in \SD_j^M(n)$, out of analogy with the $j=1, M=D^2$ case call
the embeddings $L_i : D^j \times M \to D^j \times M$ for $i\in\{1,2,\cdots,n\}$ {\it hockey pucks}. 
$L_0$ is the {\it long knot} associated with $L$ and $\sigma$ is the mapping (only well-defined modulo the equivalence
relation on splicing diagrams) from the {\it relative heights} of the pucks
to their {\it indices}, i.e. $L_i$ has height $\sigma^{-1}(i)$. $L_{\sigma(1)}$ is a bottom-most puck,
$L_{\sigma(n)}$ is a top-most puck.   
\begin{itemize}
\item[(a)] Hockey pucks allow for the construction of re-embedding maps. Given a hockey puck $L_i$ and
$f \in \EK{j,M}$, the function $L_i.f := L_i \circ f \circ L_i^{-1}$ is defined on the image of $L_i$ but we
extend the definition of $L_i.f$ to be the map $\Real^j \times M \to \Real^j \times M$ which is the identity on 
$(\Real^j\times M) \setminus L_i([-1,1]^j\times M)$.  Notice that this function can only fail to be smooth on the set $L_i(([-1,1]^j)^\circ \times \partial M)$, and generally this is precisely the set of points where $L_i.f$ fails to be differentiable.  Splicing diagrams give rise to a {\it splicing operation} (Definition \ref{def1}) and the continuity constraint allows for this to be a smooth map. 
\item[(b)] A benefit of the continuity constraint is that it makes splicing diagrams into objects that are similar to links. For example, given $(L_0,L_1,\cdots,L_n,\sigma) \in \SD_1^{D^2}(n)$, generically $(L_{0|\Real \times \{0\}}, L_{1|\{0\}\times S^1}, \cdots, L_{n|\{0\}\times S^1})$ will be a KGL. Given $(L_0,L_1,\cdots,L_n,\sigma) \in \SD_j^{D^k}(n)$, $(L_{0|\Real^j \times \{0\}}, L_{1|\{0\}\times S^{k-1}}, \cdots, L_{n|\{0\}\times S^{k-1}})$ is generically a link with one component `long'.  There are certain circumstances where these embeddings will not be disjoint. These are rare yet important cases, see Propositions \ref{cubesinclusion} and \ref{shrinkingprop}.
\item[(c)] The definition of a splicing diagram does not explicitly state that $(L_{1|\{0\}\times S^{k-1}},\- \cdots,\- L_{n|\{0\}\times S^{k-1}})$ is a trivial link when $M=D^k$, but it follows by a simple
induction argument -- by design the bottom-most hockey puck is disjoint from the other link components.  Theorem \ref{afflin} can be seen as an enhanced version of this observation.
\end{itemize}
\item[3)] For the sake of defining a single splicing operation, disjointness of the pucks is perfectly acceptable.  But there are isotopies between spliced knots (coming from diagrams with disjoint pucks) that can not be realized as splices with the pucks disjoint throughout.  By keeping track of the permutation $\sigma$ and allowing non-disjointness of pucks, the definition of splicing diagrams allows the splicing operad, as a space, to capture natural isotopies that happen in spaces of knots.  Meaning, the splicing operad more accurately reflects the homotopy-type of embedding spaces. 
\end{itemize}

\begin{eg}
An example of the action of $\SD_1^{D^2}$ on $\ED{1,D^2}$ from Definition \ref{def1}. 
{
\psfrag{f1}[tl][tl][0.8][0]{$f_1$}
\psfrag{f2}[tl][tl][0.8][0]{$f_2$}
\psfrag{p1}[tl][tl][0.8][0]{$1$}
\psfrag{m1}[tl][tl][0.8][0]{$-1$}
\psfrag{l0}[tl][tl][1][0]{$L_0$}
\psfrag{l1}[tl][tl][0.8][0]{$L_1$}
\psfrag{l2}[tl][tl][0.8][0]{$L_2$}
\psfrag{splice}[tl][tl][1][0]{$L.F = L_2.f_2 \circ L_1.f_1 \circ L_0$}
$$\includegraphics[width=12cm]{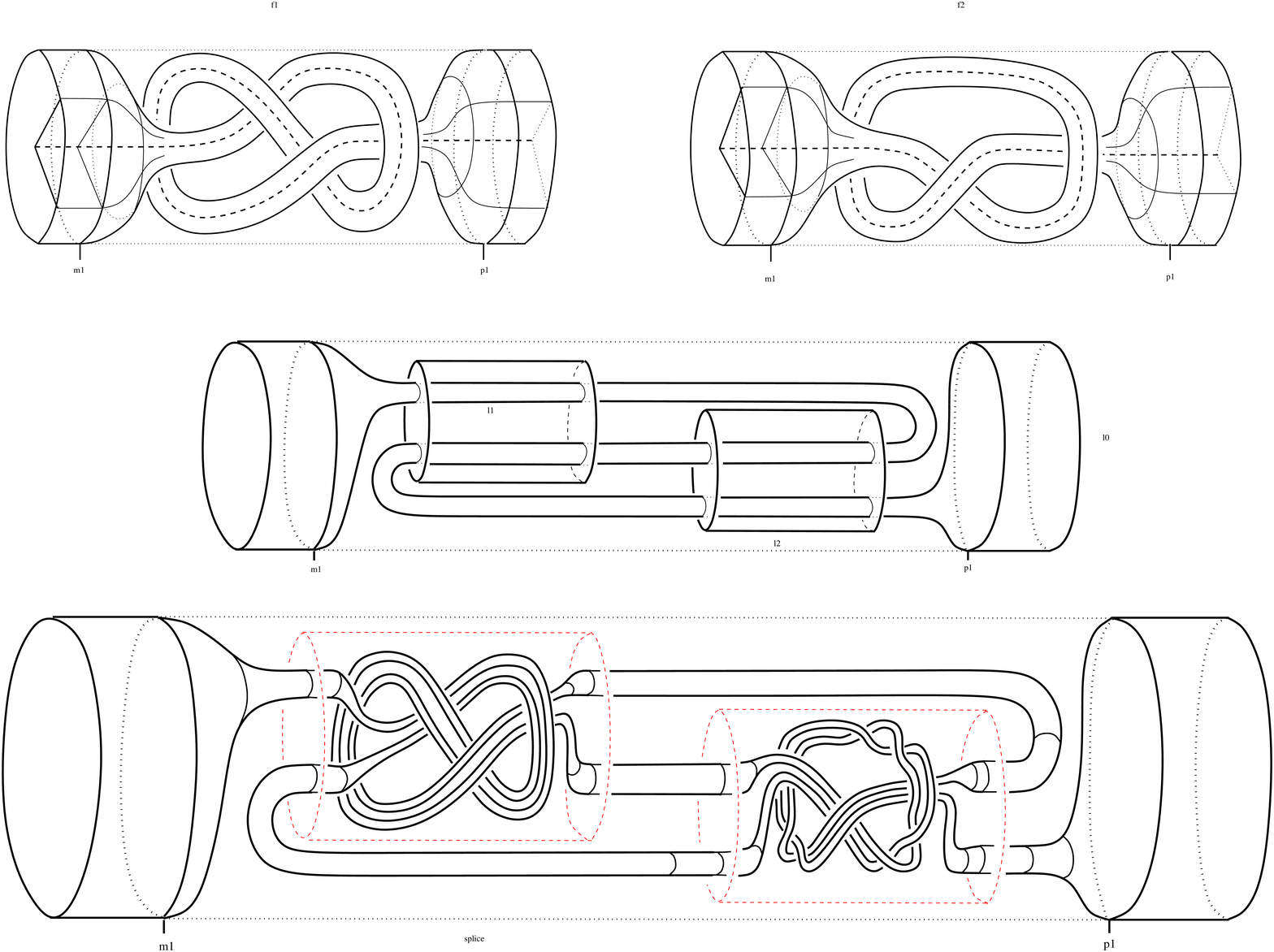}$$
}
In this example we are thinking of the figure-8 and trefoil knots as normalized to be in $\hat \K_{3,1}$, which
explains the $3$-fold twisting seen in the bottom long knot, as the trefoil's `blackboard framing' 
disagrees with its `homological framing' by three twists, while both framings are the same for the figure-8 knot.
\end{eg}

\begin{defn}\label{def1}
Let $L=(L_0,L_1,\cdots,L_n,\sigma) \in \SC_j^M(n)$ and $F=(f_1,\cdots,f_n) \in \EK{j,M}^n$. 
$$ L.F = (L_{\sigma(n)}.f_{\sigma(n)})\circ \cdots \circ (L_{\sigma(2)}.f_{\sigma(2)}) \circ (L_{\sigma(1)}.f_{\sigma(1)}) \circ L_0 \in \EK{j,M}$$
where $L_i.f_i = L_i \circ f_i \circ L_i^{-1}$ and we use the convention that $L_i.f_i$ 
is defined to be the identity outside of the image of $L_i$. $L.F$ is called the splicing operation of $L$ on $F$.
\end{defn}

The remainder of this section is devoted to showing that the space of splicing diagrams
forms an operad, and the splicing operation defined above becomes an operad action on
$\EK{j,M}$.   

Given a collection of composable functions
$$\xymatrix{ A_0 \ar[r]^{f_1} & A_1 \ar[r]^{f_2} & A_2 \ar[r]^{f_3} & \cdots \ar[r]^{f_{n-1}} & A_{n-1} \ar[r]^{f_n} & A_n }$$
their composite will be denoted 
$$\bigcirc_{i=1}^n f_i : A_0 \to A_n.$$

\begin{prop}\label{def2}
The collection $\SC_j^M = \sqcup_{n=0}^\infty \SC_j^M(n)$ is a multiplicative $\Sigma$-operad.  With Definition \ref{def1}, 
$\SC_j^M$ acts on $\EK{j,M}$.  The operad's structure map has the form 
$$\SC(k)\times \left(\SC(j_1) \times \cdots \times \SC(j_k) \right)\to \SC(j_1+\cdots+j_k)$$
(superscripts $M$ and subscripts $j$ suppressed) and is defined below. 
Let $J=(J_0,J_1,\cdots,J_k,\alpha) \in \SC(k)$ and $(L_i,\sigma_i) \in \SC(j_i)$ for $i=1,2,\cdots,k$, 
then $J.L \in \SC(j_1+\cdots+j_k)$ has $0$-th entry
$$ \left( \bigcirc_{i=1}^k (J_{\alpha(i)}L_{\alpha(i)0}J_{\alpha(i)}^{-1})\right)J_0.$$
The $(a,b)$-th coordinate entry for $a \in \{1,\cdots, k\}$ and
$b \in \{1,\cdots, j_a\}$ is given by
$$ \left( \bigcirc_{i=\alpha^{-1}(a)+1}^k (J_{\alpha(i)}L_{\alpha(i)0}J_{\alpha(i)}^{-1})\right)J_aL_{a,b}. $$
As with Definition \ref{overlappingcubes} we identify the pairs $\{(a,b) : a \in \{1,\cdots,k\}, b \in \{1,\cdots,j_a\}\}$  
with the set $\{1,\cdots,j_1+\cdots+j_k\}$ via the lexicographical ordering.
The permutation associated to $J.L$ is the natural one induced by the permutations
$(\alpha, \sigma_1, \cdots, \sigma_k)$ as in Definition \ref{overlappingcubes}.
The right action of $\Sigma_n$ on $\SC(n)$ is given by 
$$(J_0,J_1,\cdots,J_n,\alpha).\sigma = (J_0, J_{\sigma(1)}, \cdots, J_{\sigma(n)}, \sigma^{-1}\alpha).$$

\begin{proof}
(1) Associativity. For this we need to show $J.(L.M) = (J.L).M$. Let 
$M_{a,b} = (M_{a,b,0}, M_{a,b,1}, \-\cdots, M_{a,b,\beta_{a,b}}, \gamma_{a,b})$.
Notice the $(a,b,c)$-th entry of $J.(L.M)$ is given by 
$$\left( \bigcirc_{i=\alpha^{-1}(a)+1}^k J_{\alpha(i)} \left( \bigcirc_{n=1}^{j_{\alpha(i)}}
 L_{\alpha(i), \sigma_{\alpha(i)}(n)} M_{\alpha(i), \sigma_{\alpha(i)}(n), 0} L^{-1}_{\alpha(i), \sigma_{\alpha(i)}(n)}
 \right) L_{\alpha(i), 0} J^{-1}_{\alpha(i)} \right) J_a \circ $$
 $$\left( \bigcirc_{i=\sigma^{-1}(a) +1}^{j_a} 
 L_{a, \sigma_a(i)} M_{a, \sigma_a(i), 0} L^{-1}_{a, \sigma_a(i)} \right) L_{a,b} M_{a,b,c}$$
 while the $(a,b,c)$-th entry of $(J.L).M$ is given by

$$ \left( \bigcirc_{(i,n)>(\alpha^{-1}(a), \sigma^{-1}_a(b))} \left( \bigcirc_{l=i+1}^k J_{\alpha(l)}L_{\alpha(l), 0} J_{\alpha(l)}^{-1} \right) J_{\alpha(i)}L_{\alpha(i), \sigma_{\alpha(i)}(n)} M_{\alpha(i), \sigma_{\alpha(i)}(n), 0}
 L_{\alpha(i), \sigma_{\alpha(i)}(n)}^{-1} J_{\alpha(i)}^{-1} \circ \right.$$
$$\left. \left( \bigcirc_{l=k}^{i+1} J_{\alpha(l)}L_{\alpha(l), 0}^{-1} J_{\alpha(l)}^{-1} \right)\right) \left( \bigcirc_{i=\alpha^{-1}(a)+1}^k J_{\alpha(i)}L_{\alpha(i), 0}J_{\alpha(i)}^{-1} \right) J_a L_{a,b} M_{a,b,c}.$$

In this latter composite there are many occurrences of adjacent maps that are the inverses of each other. Cancelling these maps we see the above two expressions for the $(a,b,c)$-th term
of $(J.L).M$ and $J.(L.M)$ are identical.  Showing the $0$-th entries agree is similar. 

(2) Symmetry/Equivariance. There are two types, the `internal' equivarance, and the `external' one.  For the internal equivariance, we need to show that if $J \in \SC(k)$ and if
$L_i \in \SC(j_i)$ for all $i \in \{1,2,\cdots, k\}$ with $L=(L_1,\cdots,L_k)$ then
whenever $\sigma \in \Sigma_k$ $(J.\sigma).L = (J.(\sigma.L)).\overline{\sigma}$ where
$\sigma.L = (L_{\sigma^{-1}(1)}, \cdots, L_{\sigma^{-1}(k)})$, and $\overline{\sigma} \in \Sigma_{j_1+\cdots+j_k}$ is the associated block permutation.  This is immediate.  For the external equivarance, we need to show that if $\tau \in \Sigma_{j_1} \times \cdots \times \Sigma_{j_k}$, $J.(L.\tau)=(J.L).\tau$, which is also immediate. 

(3) Identity/Unit. The identity element in $I \in \SC(1)$ is $(Id_{\Real^j\times M}, Id_{I^j\times M}, e)$ where $e \in \Sigma_1$ is the identity element.  Given $L \in \SC(j)$ the identity axiom requires $I.L = L$ and $L.(I,I, \cdots, I) = L$, which are both satisfied.

That Definition \ref{splicedef} gives an action of $\SC$ on $\EK{j,M}$ is a special case of the above arguments, 
since the structure maps for $\SC$,  
$$\SC(k) \times \left(\SC(0) \times \cdots \times \SC(0)\right) \to \SC(0)$$
{\it is} the action of $\SC_j^M$ on $\EK{j,M}$, as $\SC_j^M(0) = \EK{j,M}$. 

A multiplicative operad is one that contains the associative operad as a sub-operad.  For 
$\SC_j^M$, the suboperad is $\{ (Id_{\Real^j \times M}, Id_{[-1,1]^j\times M}, \cdots, Id_{[-1,1]^j \times M}, Id_{\{1,2,\cdots,k\}}) : k \in \Nat \}\subset \SC_j^M$.
\end{proof}
\end{prop}

\begin{eg} An example of the structure map of $\SD_1^{D^2}$, in pictures.
{
\psfrag{l10}[tl][tl][0.8][0]{$L_{1,0}$}
\psfrag{l11}[tl][tl][0.8][0]{$L_{1,1}$}
\psfrag{l12}[tl][tl][0.8][0]{$L_{1,2}$}
\psfrag{l20}[tl][tl][0.8][0]{$L_{2,0}$}
\psfrag{l21}[tl][tl][0.8][0]{$L_{2,1}$}
\psfrag{j0}[tl][tl][0.8][0]{$J_{0}$}
\psfrag{j1}[tl][tl][0.8][0]{$J_{1}$}
\psfrag{j2}[tl][tl][0.8][0]{$J_{2}$}
\psfrag{p1}[tl][tl][0.8][0]{$1$}
\psfrag{m1}[tl][tl][0.8][0]{$-1$}
\psfrag{splice}[tl][tl][1][0]{$J.L$}
$$\includegraphics[width=12cm]{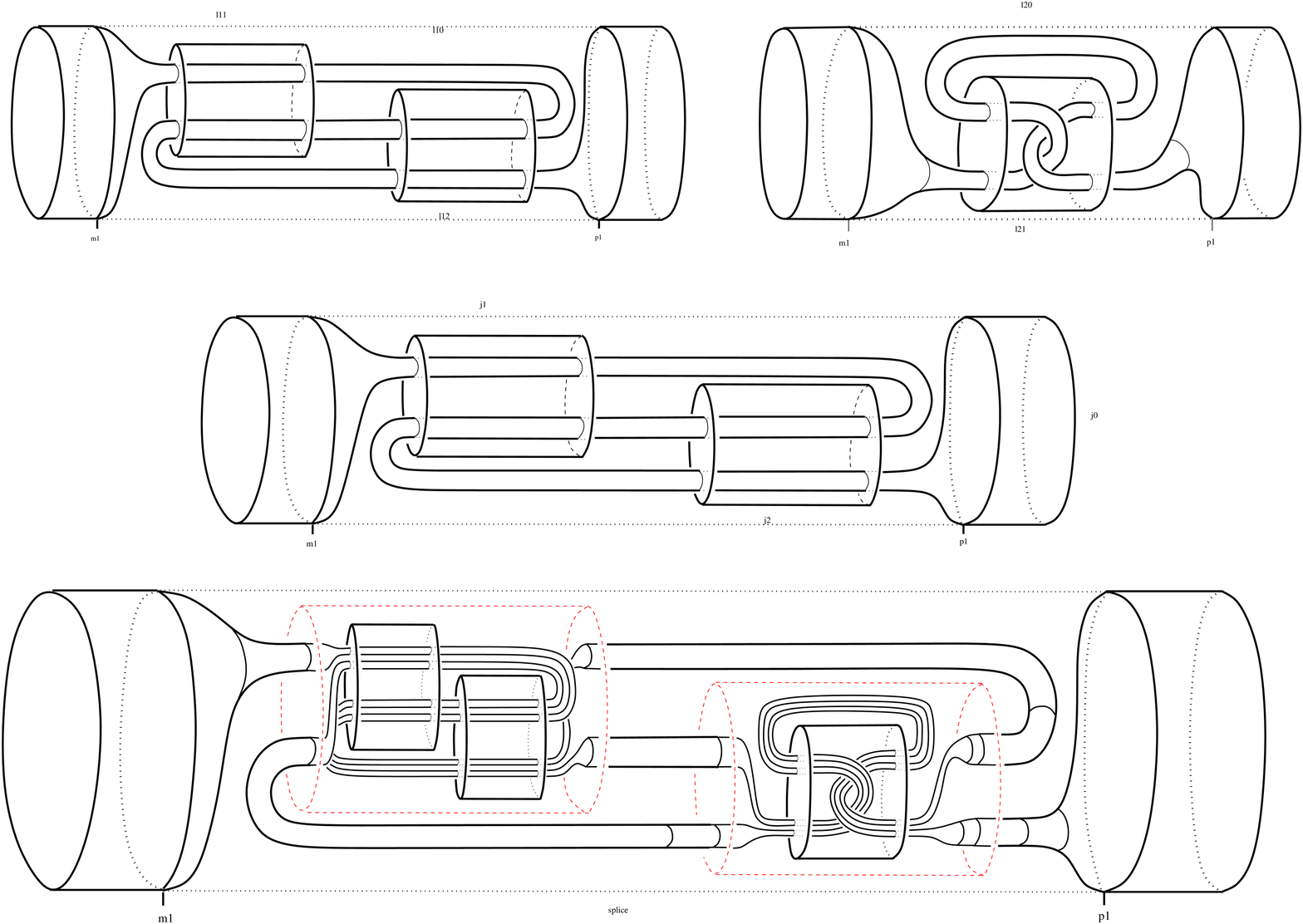}$$
}
\end{eg}

\begin{prop}\label{cubesinclusion}
There is an inclusion of operads 
$$\Cu_j' \to \SC_j^M$$
given by the maps $\Cu_j'(k) \to \SC_j^M(k)$ which have the form
$(L_1,\cdots,L_k,\sigma) \longmapsto (L_1 \times Id_M, \cdots, L_k \times Id_M, \sigma)$.
Moreover, the action of $\SC_j^M$ on $\EK{j,M}$ restricts to the action of 
Proposition \ref{flataction}.
\end{prop}

\begin{eg} The inclusion $\Cu'_1(3) \to \SC_1^{D^2}(3)$ in a picture.
{
\psfrag{p1}[tl][tl][0.8][0]{$1$}
\psfrag{m1}[tl][tl][0.8][0]{$-1$}
\psfrag{splice}[tl][tl][1][0]{}
$$\includegraphics[width=8cm]{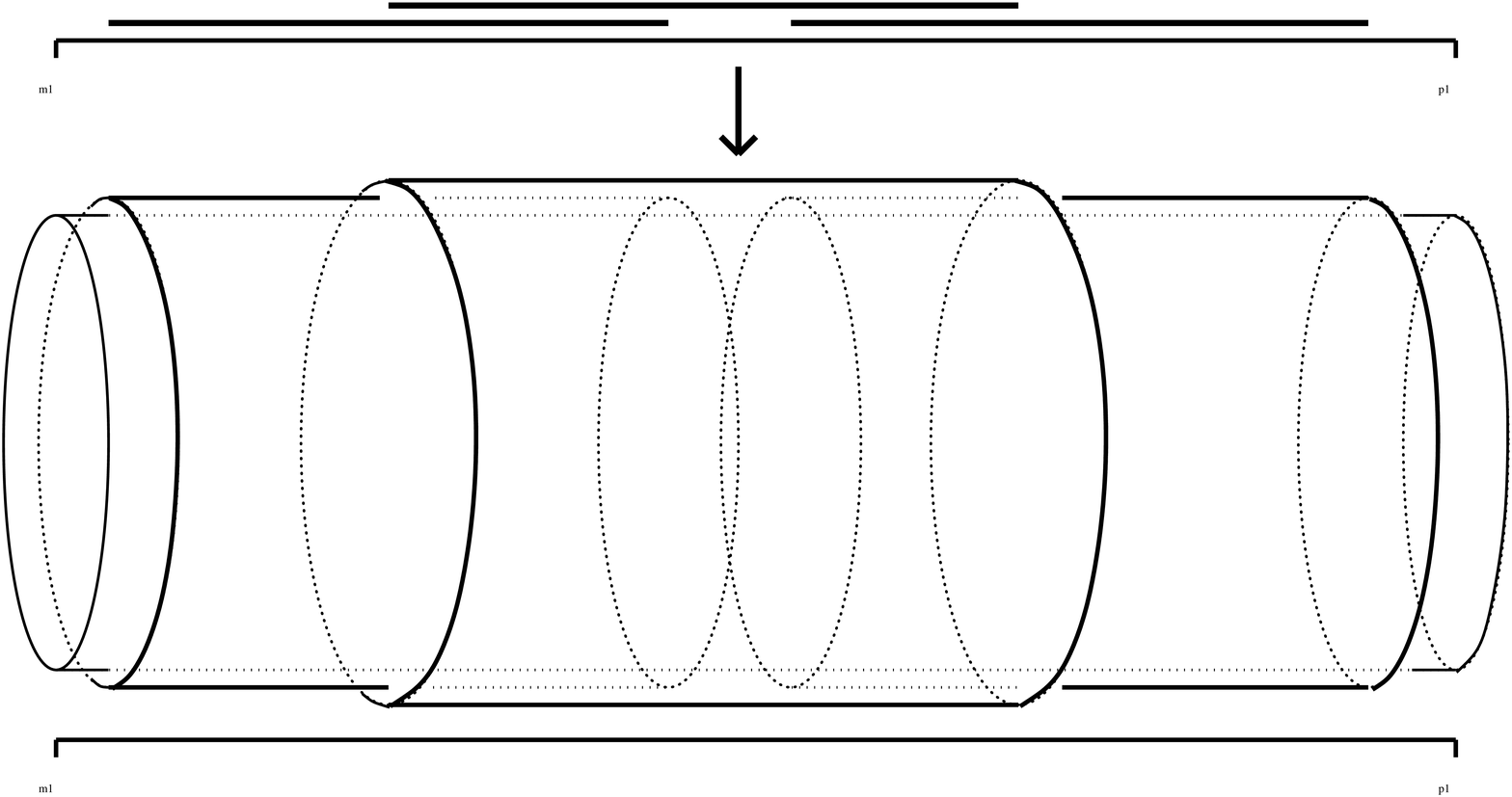}$$
}
We visualize the overlapping nature of the intervals as an infinitesimal separation in an orthogonal 
direction.  Similarly for elements of $\SC_1^{D^2}$, although we have run out of extra dimensions, so 
we depict the relative order as if one cylinder were a thin film over the other(s). 
\end{eg}

It is appealing to think of the operad $\SC_j^M(k)$ as an enhanced space of $(k+1)$-component
links where the $0$-th component is `long.'  The next proposition makes this a little more concrete in the
case that $M$ is connected with non-empty boundary.   

\begin{prop}\label{shrinkingprop}
Let $M$ be a compact connected manifold with $\partial M$ non-empty. 
$$(\SC_j^M)^\circ(k) = \{ (L_0, L_1, \cdots, L_k, \sigma) \in \SC_j^M(k) : 
(L_{0| \Real^j \times M^\circ}, L_{1| I^j \times \partial M}, \cdots, L_{k | I^j \times \partial M})
\text{ are disjoint} \}$$
Then $(\SC_j^M)^\circ = \sqcup_k (\SC_j^M)^\circ(k)$ is a suboperad without identity of $\SC_j^M$, moreover
the inclusion $(\SC_j^M)^\circ \to \SC_j^M$ is a homotopy-equivalence.
\begin{proof}
The proof is by constructing a homotopy-inverse of the inclusion $(\SC_j^M)^\circ \to \SC_j^M$.  Since it will be useful in Theorem \ref{afflin} we develop the case $M = D^n$ explicitly. Let $\beta : \Real \to \Real$ be a $C^\infty$-smooth function such that $\beta(0)=0$, $\beta'(t) \geq 0$ for all $t \geq 0$, $\beta(t)=1$ for all $t \geq 1$, $\beta(-t)=\beta(t)$ for all $t$ and $0 < \beta(t) < 1$ for $0 < t < 1$. The {\it standard shrinking map} for $\Real^j \times D^n$ is the family $ R : [0,1] \times \Real^j \times D^n \to \Real^j \times D^n $ given by $ R(t, x, v) = (x, \left( t+(1-t)\beta(|x|^2)\right) v )$. We let $R_t : \Real^j \times D^n \to \Real^j \times D^n$ denote $R(t, \cdot)$. Notice that when $t \in (0,1]$, $R_t \in \EK{j,D^n}$. 

Given $L \in \SC_j^{D^n}(k)$ and $t \in (0,1]$, let $R_t.L \in \SC_j^{D^n}(k)$ denote the
element where the $0$-th entry has the form
$$ (x,v) \longmapsto \left( \bigcirc_{i=1}^k L_{\sigma(i)}R_{t^2}L_{\sigma(i)}^{-1}\right) \circ L_0 \circ R_t( x, v ) $$
the $a$-th element has the form
$$ (x,v) \longmapsto \left( \bigcirc_{i=\sigma^{-1}(a)+1}^k L_{\sigma(i)}R_{t^2}L_{\sigma(i)}^{-1}\right) \circ 
L_a \circ R_t( tx, v)$$ $R_{1/2} : \SC_j^{D^n} \to (\SC_j^{D^n})^\circ$ is our desired homotopy-inverse. 
The general case $M \neq D^n$ proceeds similarly, using a collar neighbourhood of $\partial M \subset M$ as a
replacement for the linear structure on $D^n$. 
\end{proof}
\end{prop}

Notice that the identity element of $\SC_1^{D^2}$ is mapped via $R_{1/2}$ to the `Hopf link' in $(\SC_1^{D^2})^\circ(1)$.

The discovery of the operad $\SC^M_j$ came about fairly naturally.   Individual splicing diagrams analogous to elements of $\SC^{D^2}_j$ first appear in \cite{BudJSJ} as a formally convenient way to encode splicing.  As a topological space something similar to $\SC^{D^2}_j$ appears in \cite{KS} when describing the homotopy-type of various components of
$\K_{3,1}$. Thus ideas similar to Definition \ref{def1} have been present for some time.  Now consider making Definition \ref{def1} satisfy an associativity law for a hypothetical operad structure on $\SC_j^M$.   Since the associativity law for an operad action uses the structure map of an operad {\it only once} and the action of the operad on $\EK{j,M}$ {\it three times}, one could use the associativity condition together with a hypothetical action in an attempt to intuit an operad structure map $\SC_j^M(k) \times \prod_{i=1}^k\SC_j^M(j_i) \to \SC_j^M(j_1+\cdots+j_k)$. This works and is precisely how the author was led to define the operad structure maps for $\SC^M_j$. 

\begin{defn}
We denote the wreath product of a group $G$ and $\Sigma_n$ be $\Sigma_n \wr G$. The main purpose of the wreath product for this paper is that it is the appropriate group that extends two natural group actions.  If $G$ acts on $X$, $G^n$ acts on $X^n$ via the product action and $\Sigma_n$ acts on $X^n$ via the regular representation. $\Sigma_n \wr G$ fits into a short exact sequence $0 \to G^n \to \Sigma_n \wr G \to \Sigma_n \to 0$.  Moreover, $\Sigma_n \wr G$ acts on $X^n$ and its action is equivariant with respect to this short exact sequence.  $\Sigma_n \wr G$ is the group $\Sigma_n \ltimes G^n = \mathrm{Aut}\{1,2,\cdots,n\} \ltimes G^{\{1,2,\cdots,n\}}$, i.e. the semi-direct product of $G^n$ and $\Sigma_n$ where $\Sigma_n$ acts on $G^n$ by the regular representation.  We will use the notation $\Sigma_n^* \wr G$ to denote $G \times (\Sigma_n \wr G)$. $\Sigma^*_n \wr G$ should be thought of as the above wreath product construction but with the identification $\Sigma_n^* = \mathrm{Aut}(\{0,1,2,\cdots,n\} \text{ fixing } 0)$, i.e. $\Sigma_n^* \wr G = \mathrm{Aut}(\{0,1,\cdots,n\} \text{ fixing } 0) \ltimes G^{\{0,1,\cdots,n\}}$. We denote the sequence of groups $\sqcup_n \Sigma_n^* \wr G$ by $\Sigma^* \wr G$.  Since a preferred copy of $G$ splits off $\Sigma^* \wr G$, if $X$ is a space with an action of $\Sigma^* \wr G$, when $g \in G$ and $x \in X$, $g.x \in X$ will refer to the action of $G$ on $X$ coming from this preferred factor. 

A $\Sigma^* \wr G$-operad $\mathcal O$ is a sequence of spaces $\mathcal O(n)$ for $n \in \Nat$ together with group actions of $\Sigma^*_n \wr G$ on $\mathcal O(n)$ for all $n \in \Nat$ satisfying an (1) associativity axiom, a (2) symmetry axiom and an (3) identity axiom.  The (1) associativity and (3) identity axioms are exactly as in the definition of a $\Sigma$-operad.  The symmetry axiom (2) has two parts, an `inner' equivariance, together with an `outer' equivariance.  

The `inner' equivariance condition can be expressed as $(J.\gamma).L = (J.(\overline{\gamma}.L)).\tilde \gamma$, where
$\gamma \in \Sigma^*_k \wr G$, $J \in \mathcal O(k)$ and $L \in \prod_{i=1}^k \mathcal O(j_i)$.  If we write $\gamma = (g_0, \beta, g_1, \cdots, g_k)$ and $L \in \prod_{i=1}^k \mathcal O(j_i)$ as $(L_1, \cdots, L_k)$ then 
$$\overline{\gamma}.L = (g_{\beta^{-1}(1)}.L_{\beta^{-1}(1)}, \cdots, g_{\beta^{-1}(k)}.L_{\beta^{-1}(k)}).$$ 
Similarly, if $H \in \mathcal O(j_1 + \cdots + j_k)$ then $H.\tilde \gamma = g_0^{-1}.H.\overline{\beta}$
where $\overline{\beta} \in \Sigma_{j_1 + \cdots + j_k}$ is the block permutation associated to $\beta \in \Sigma_k$. 

The outer equivariance can be expressed as $J.(L.\gamma) = (J.L).\gamma$ whenever $\gamma = \gamma_1 \times \cdots \times \gamma_k$, $\gamma_i \in \Sigma_{j_i} \wr G$.  Note we do not allow $\gamma_i \in \Sigma^*_{j_i} \wr G$ for this condition. 
\end{defn}

The next proposition will investigate further equivariance properties of the splicing operads and their actions. Let $Diff(I^j \times M)$ and denote the group of diffeomorphisms of $I^j \times M$ that restrict to diffeomorphisms of $(\partial I^j)\times M$, where $I = [-1,1]$. 
Similarly, let $Diff(D^j \times M)$ be the diffeomorphisms of $D^j \times M$ that restrict to diffeomorphisms of $(\partial D^j) \times M$. 

\begin{prop}\label{equivarianceprop}
$\EK{j,M}$ is taken to be a $Diff(I^j \times M)$-space, where the action is by conjugation.  Similarly, $\ED{j,M}$ is a $Diff(D^j \times M)$-space.  There is an action of $\Sigma^* \wr Diff(I^j \times M)$ on $\SC_j^M$ making $\SC_j^M$ into a $\Sigma^* \wr Diff(I^j \times M)$-operad.  Similarly, there is an action of $\Sigma^* \wr Diff(D^j \times M)$ on $\SD_j^M$ making $\SC_j^M$ into a $\Sigma^* \wr Diff(D^j \times M)$-operad.  $\SC_j^M$ and $\SD_j^M$ act on $\EK{j,M}$ and $\ED{j,M}$ in the sense of $\Sigma^* \wr G$-operad actions. 

\begin{proof}
The right action of $\Sigma^*_k \wr Diff(D^j\times M)$ on $\SD_j^M(k)$ is given by:
$$\SD_j^M(k) \times \Diff(D^j \times M )\times \left(\Sigma_k \ltimes \Diff(D^j \times M)^k\right) \to \SD_j^M(k)$$
$$(J_0, J_1,\cdots,J_k, \sigma), g_0, (\gamma, g_1, \cdots, g_k) \longmapsto 
(g_0^{-1} \circ J_0 \circ g_0, g_0^{-1} \circ J_{\gamma(1)} \circ g_{1}, \cdots, 
                               g_0^{-1} \circ J_{\gamma(k)} \circ g_{k}, \gamma^{-1}\sigma)$$
Abbreviate $J = (J_0, J_1,\cdots,J_k, \sigma)$ and $g = ( g_0, \gamma, g_1, \cdots, g_k)$.  Let 
$L = (L_1, \cdots, L_k) \in \prod_{i=1}^k \SD_j^M(j_i)$, and write
$L_i = (L_{i \ 0}, L_{i \ 1}, \cdots, L_{i \ j_i}, \alpha_i)$.  Then $(J.g).L \in \SD_j^M(\sum_i j_i)$, whose
$0$-th entry is 
$$\left( \bigcirc_{i=1}^k g_0^{-1}J_{\sigma(i)}g_{\beta^{-1}\sigma(i)}L_{\beta^{-1}\sigma(i) \ 0}g_{\beta^{-1}\sigma(i)}^{-1}J_{\sigma(i)}^{-1}g_0\right) g_0^{-1} J_0 g_0$$ 
and whose $(a,b)$-th entry (before lexicographically ordering) is
$$\left( \bigcirc_{i=\sigma^{-1}\beta(a)+1}^k g_0^{-1}J_{\sigma(i)}g_{\beta^{-1}\sigma(i)}L_{\beta^{-1}\sigma(i) \ 0}g_{\beta^{-1}\sigma(i)}^{-1}J_{\sigma(i)}^{-1}g_0\right) g_0^{-1} J_{\beta(a)}g_a L_{a,b}$$
cancelling inverse maps, these two expressions reduce to 
$$g_0^{-1}\left( \bigcirc_{i=1}^k J_{\sigma(i)} \left(g_{\beta^{-1}\sigma(i)}L_{\beta^{-1}\sigma(i) \ 0}g_{\beta^{-1}\sigma(i)}^{-1}\right) J_{\sigma(i)}^{-1}\right) J_0 g_0$$ 
and 
$$g_0^{-1}\left( \bigcirc_{i=\sigma^{-1}\beta(a)+1}^k J_{\sigma(i)} \left(g_{\beta^{-1}\sigma(i)}L_{\beta^{-1}\sigma(i) \ 0}g_{\beta^{-1}\sigma(i)}^{-1}\right)J_{\sigma(i)}^{-1}\right) J_{\beta(a)}g_a L_{a,b}$$
respectively, which are the entries of $(J.(\overline{g}.L)).\tilde g$. 
The `outer' equivariance condition is immediate. 
\end{proof}
\end{prop}
 
\section{The homotopy type of the splicing operad}

The next theorem should be thought of as a semi-linear ordering enhancement of Cerf's 
homotopy-classification of spaces of tubular neighbourhoods \cite{Cerf}.  

\begin{thm}\label{afflin}
Let $\mathcal{LO}_{j,n}(k) \subset \SD_j^{D^n}(k)$ be the subspace where the embeddings
$L_i : D^j \times D^n \to D^j \times D^n$ are affine linear for $i \in \{ 1,2, \cdots, k \}$. 
Then the inclusion $\mathcal{LO}_{j,n}(k) \to \SD_j^{D^n}(k)$ is a homotopy-equivalence 
for all $k \in\{1,2,3\cdots\}$. 
\begin{proof}
Recall the standard shrinking map from the proof of Proposition \ref{shrinkingprop}. 
Given $L \in \SD_j^{D^n}(k)$ and $t \in (0,1]$, let $R_t.L \in \SD_j^{D^n}(k)$ denote the
element where the $0$-th entry has the form
$$ (x,v) \longmapsto \left( \bigcirc_{i=1}^k L_{\sigma(i)}R_{t^2}L_{\sigma(i)}^{-1}\right) \circ L_0 \circ R_t( x, v ) $$
the $a$-th element has the form
$$ (x,v) \longmapsto \left( \bigcirc_{i=\sigma^{-1}(a)+1}^k L_{\sigma(i)}R_{t^2}L_{\sigma(i)}^{-1}\right) \circ 
L_a \circ R_t( tx, v)$$

The idea of the proof is to shrink elements $L \in \SD_j^{D^n}(k)$ to the point where we can apply a linearization process.  The linearization process $[0,1] \times D^j \times D^n \to \Real^{j+n}$ 
applied to $L_i$ for $i \in \{1,2,\cdots, k\}$ is given by
$$(t,x,v) \longmapsto 
\left\{
\begin{array}{ll}
\frac{1}{t}\left(L_i( t(x,v) ) - L_i(0,0)\right) + L_i(0,0) & 0 < t \leq 1\\
(DL_i)_{(0,0)}(x,v) + L_i(0,0)              & t=0
\end{array}
\right..
$$

If we think of this as a time-varying family of maps $L_{i t} : D^j \times D^n \to \Real^{j+n}$, we can make some observations on the family.
\begin{itemize}
\item[(a)] For all $t$ the map $L_{it} : D^j \times D^n \to \Real^{j+n}$ is an embedding, thus the family is an isotopy of $L_i$. 
\item[(b)] $L_i$ and $L_{i t}$ are uniformly close, moreover, an upper bound on their $C^0$-distance is given by the maximum of the norm of the Hessian of $L_i$.  
\item[(c)] Under the shrinking map the 2nd derivative of $L_i$ goes to zero at an order of magnitude faster than the 1st derivative. 
\end{itemize}

Given any $L \in \SD_j^{D^n}(k)$, we can apply the shrinking map until linearization can be applied to the $(L_1, \cdots, L_k)$ part of the family.  Via linearization we can ensure $(L_{1 t | D^j \times \partial D^n}, \cdots, L_{k t | D^j \times \partial D^n})$ are disjoint.  Apply isotopy extension to the isotopy $(L_{1 t | D^j \times \partial D^n}, \cdots, L_{k t | D^j \times \partial D^n})$ allows us to construct the family $L_{0 t}$.  This gives us a path in $\SD_j^{D^n}(k)$ that begins
at $L$ and ends in $\mathcal{LO}_{j,n}(k)$.  Moreover, we choose how long to run the shrinking map based on the maximum of the 2nd derivative of $L$, which varies continuously on $\SD_j^{D^n}(k)$. Similarly the isotopy extension, since it is a solution to an ODE varies continuously with the input isotopy.  This gives us a homotopy of the identity map on $\SD_j^{D^n}(k)$ to a map $\SD_j^{D^n}(k) \to \mathcal{LO}_{j,n}(k)$, which is a homotopy-inverse to the inclusion $\mathcal{LO}_{j,n}(k) \subset \SD_j^{D^n}(k)$.  
\end{proof}
\end{thm}

There is a related theorem of Brendle and Hatcher \cite{BH}, who have shown that in dimension $3$ the space of unlinks has the homotopy-type of the subspace of round unlinks.  Their proof is analogous, one key difference is their step where they add spanning discs to their trivial links -- this is via an application of the theorem that $Diff(S^3) \simeq O_4$.  To make the analogy a little more explicit, the shrinking construction above supplies a homotopy-equivalence between $\SD_j^{D^n}(k)$ and a subspace of $\SD_j^{D^n}(k)$ where each $L$ has a unique semi-linear ordering $\sigma \in \Sigma_k$ up to equivalence (this is essentially the `separated' subspace in \cite{BH}).  This subspace of $\SD_j^{D^n}(k)$ is therefore a genuine embedding space and therefore has the homotopy-type of a CW-complex \cite{HW}. 

\section{Splicing classical knots}

The point of this section is to show how the splicing operad is in some sense a more natural
operad than cubes operads for the purposes of describing the homotopy-type of embedding spaces.   This is largely
done by example, for the splicing operad's action on the space $\K_{3,1}$.  We start by refining  
the splicing operad $\SD_1^{D^2}$, throwing away the parts that contain redundant information from the
point of view of the action on $\K_{3,1}$, to produce the irreducible splicing operad $\SP_{3,1}$. 
We then show $\K_{3,1}$ to be free over $\SP_{3,1}$. Further we show $\SP_{3,1}$ to be a free product 
of $\Cu_1' \rtimes O_2$ and a free operad over a $\Sigma^* \wr O_2$-space, which we identify.

Let $\hat \K_{3,1} \subset \ED{1,D^2}$ be the subspace with zero homological framing from Definition 
\ref{longdef}. $\SD_1^{D^2}$ acts on $\ED{1,D^2}$ but notice that it does {\it not}
restrict to an action on $\hat \K_{3,1}$ since it does not preserve the homological framing
of the knot.  Moreover, not every element of $\SD_1^{D^2}$ results in a useful splicing construction -- think for example of an element $(L_0,L_1,\sigma) \in \SD_1^{D^2}(1)$ where $L_1$ is disjoint from $L_0$.
Below we define a suitable suboperad of $\SD_1^{D^2}$ that acts on $\hat \K_{3,1}$ in 
a useful way. 

\begin{defn}
The $3$-dimensional irreducible splicing operad $\SP_{3,1}$ is the subset of $\SD_1^{D^2}$ where
$(L,\sigma) \in \SD_1^{D^2}(k)$ is an element of $\SP_{3,1}(k)$ provided all of the following conditions are satisfied:
\begin{itemize}
\item[1)] $\SP_{3,1}(0) = \emptyset$, i.e. this is an operad with empty base.
\item[2)] We demand that $L_i$ is an orientation-preserving embedding for each $i\in \{1,2,\cdots,k\}$.
\item[3)] $L_0 \in \hat \K_{3,1}$, meaning that the linking numbers of
$L_{0|\Real \times \{(0,0)\}}$ and $L_{0|\Real \times \{(1,0)\}}$ are zero. 
\item[4)] The link corresponding to $L$ is irreducible. 
\item[5)] Every incompressible torus in the complement of the
link associated to $L$ separates components of $L$. 
\end{itemize}
\end{defn}

Condition (4) above uses {\it irreducible} in the sense of knot theory, that one can not separate components
of the link $$(L_{0 | \Real^j \times \{0\}}, L_{1|\{0\}\times S^{n-j-1}}, \cdots, L_{k|\{0\}\times S^{n-j-1}})$$ 
by embedded co-dimension zero balls. It can be restated as saying that the path-component of 
$(L,\sigma)$ in $\SD_j^{D^{n-j}}(k)$ does not contain a representative $(L',\sigma')$ such that for some $i \in\{1,2,\cdots,k\}$ $L'_i$ is disjoint from $L'_0$.  Conditions (1) and (5) can be restated as saying the JSJ-decomposition of the complement of $L$ contains no knot complements (only link complements with two or more components are permitted in the JSJ-decomposition).  Note also that condition (4) forces condition (1), since if the base of the operad were non-empty, the resulting degeneracy maps (see the comments following Definition \ref{operadDEF}) could produce reducible links, as in the case of the Borromean rings thought of as an element of $\SP_{3,1}(2)$. 

It's interesting to consider how one might want to generalize the irreducible splicing operad to an appropriate
irreducible splicing operad $\SP_{n,j}\subset \SD_j^{D^n}$ for all $n$ and $j \geq 1$.  There appears to be no high-dimensional analogue of (3).  Condition (4) immediately generalizes, although it's not clear when splicing preserves (4).  The natural generalization of (5) would be to talk about incompressible $S^j \times S^{n-j-1}$ manifolds in the link complement, presumably where incompressible means not bounding a $D^{j+1} \times S^{n-j-1}$, although perhaps a more flexible definition would be desireable. 

By the work of Hatcher \cite{Hatcher2}, $Diff(D^1 \times D^2)$ has the homotopy-type
of its linear subgroup $O_2\times \Zed_2$. The subgroup that preserves the orientation of
$D^1 \times D^2$ is isomorphic to $O_2$, so we can consider $\SP_{3,1}$ to be a $\Sigma^*\wr O_2$-operad
and $\hat \K_{3,1}$ as a space with an $O_2$-action given by conjugation. Note that with this action action of $O_2$ on $D^1 \times D^2$, mirror reflections reverse the orientations of both $D^1$ and $D^2$ factors. 

\begin{defn}\label{complexitydef}
Given $(L,\sigma) \in \SP_{3,1}(k)$, let $\hat L \subset S^3$ denote the associated link in $S^3$. 
The idea is to consider $S^n$ as the one-point compactification of $\Real^n$. $\hat L$ has $(k+1)$-components
$\hat L_0$ is the one-point compactification of $L_{0|\Real\times\{0\}} : \Real \to \Real^3$. 
$\hat L_i$ is the image of $L_{i|\{0\}\times S^1} : S^1 \to [-1,1]\times D^2 \subset \Real^3 \subset S^3$. 
Given $(L,\sigma) \in \SP_{3,1}(k)$ we say it is Seifert or hyperbolic respectively
if the associated link $\hat L \subset S^3$ has Seifert-fibred or hyperbolic complement, respectively. 

Given a $3$-manifold $M$ let $c(M)$ denote the number of components of $M$ split
along its canonical (geometric) decomposition.  We ignore the compression-body decomposition. 
So for a knot $K$ in $S^3$, the complexity of its complement $c(K)$ is $0$ if and only if it is the unknot (since after compression the manifold is empty), 
$1$ if and only if it is a torus or hyperbolic knot.  Similarly for $L$ an irreducible KGL, 
$c(L) = 0$ if and only if $L$ is the unknot, $c(L) = 1$ if and only if $L$ is hyperbolic or 
Seifert. 

Given a link $L$ in $S^3$, the symmetry group of the link is denoted $\pi_0 Diff(S^3,L)$, i.e.
the mapping class group of the pair $(S^3,L)$. Given $L \in \SP_{3,1}(k)$, the symmetry group
$B_L$ of $L$ is the defined to be a subgroup of $\pi_0 Diff(S^3, \hat L)$, 
where we put the additional restriction that the action on $S^3$ is by orientation-preserving 
diffeomorphisms and we require that the $\hat L_0$ component is preserved.  
\end{defn}

\begin{prop} \cite{BudJSJ}
The splicing map 
$$\SP_{3,1}(k) \times \prod_{i=1}^k \SP_{3,1}(j_i) \to \SP_{3,1}(j_1+ \cdots + j_k)$$ 
satisfies
$$c(J.(L_1, \cdots, L_k)) = c(J) + \sum_{i=1}^k c(L_i)$$
except in the two possible degenerate cases:
\begin{itemize}
\item[(a)] $\hat J$ is a Hopf link, or $\hat L_i$ is a Hopf link for some $i$.
\item[(b)] $\hat J$ contains two parallel components, $\hat J_a$ and $\hat J_b$ for $a,b>0$, i.e. 
$\hat J_a$ and $\hat J_b$ bound an untwisted embedded annulus disjoint from 
$(\hat J_0 \cup \hat J_1 \cup \cdots \cup \hat J_k) \setminus (\hat J_a \cup \hat J_b)$, and
either $\hat L_a$ or $\hat L_b$ are not prime with respect to connect-sum along the $0$-th strand. 
\end{itemize}
\end{prop}

If $c(J.(L_1, \cdots, L_k)) = c(J) + \sum_{i=1}^k c(L_i)$ we call $J.L$ a {\it non-redundant splice}.
$c(J.(L_1, \cdots, L_k)) < c(J) + \sum_{i=1}^k c(L_i)$ type redundant splices are the only ones
possible in the splicing operad $\SP_{3,1}$ \cite{BudJSJ}.  In the larger operad $\SD_1^{D^2}(k)$
redundant splices of the form $c(J.(L_1, \cdots, L_k)) > c(J) + \sum_{i=1}^k c(L_i)$ are possible,
but this requires one of $\{L_1,\cdots,L_k\}$ to be the unknot. 

Every (isotopy class of) element of $\hat\K_{3,1}$ and $\SP_{3,1}$ can be expressed as an iterated non-redundant 
splice of objects from $\hat\K_{3,1}$ and $\SP_{3,1}$ whose complements $M$ satisfy $c(M)=1$. 
Moreover, up to isotopy and the action of $\Sigma^* \wr O_2$ on $\SP_{3,1}$, this decomposition is unique \cite{BudJSJ}. 
This should be thought of as the analogous unique decomposition theorem to Schubert's prime factorization of
knots, but for satellite operations.  Theorems \ref{freenessthm} and \ref{operadstruc} give the generalization of the above to a statement about the homotopy-type of spaces of knots. 

\begin{thm}\label{freenessthm}
Let $\mathcal{TH} \subset \hat\K_{3,1}$ be the subspace consisting of knots which are either non-trivial torus knots, or hyperbolic knots. Then the action of $\SP_{3,1}$ on $\hat \K_{3,1}$ induces an $O_2$-equivariant homotopy-equivalence 
$$ \SP_{3,1}(\mathcal{TH}) \equiv \sqcup_{j=0}^\infty \left(\SP_{3,1}(j) \times_{\Sigma_j \wr O_2} \mathcal{TH}^j\right) \to \hat\K_{3,1}.$$
The action of $O_2$ on $\SP_{3,1}(\mathcal{TH})$ is induced by the outer action of $O_2$ on $\SP_{3,1}(j)$.  The action of $\Sigma_j \wr O_2$ on $\SP_{3,1}(j)$ is given by the inner action (see Proposition \ref{equivarianceprop}). Further, the components of $\mathcal{TH}$ have two possible homotopy-types:
\begin{itemize}
\item[(a)] A torus knot component of $\mathcal{TH}$ has the homotopy-type of $S^1$. If $f \in \hat \K_{3,1}$ is a torus knot there is an $O_2$-equivariant homotopy-equivalence $S^1 \to \hat \K_{3,1}(f)$. The action of $O_2$ on $S^1$ is standard. $\hat \K_{3,1}(f)$ denotes the path-component of $\hat \K_{3,1}$ containing $f$.
\item[(b)] A hyperbolic knot component of $\mathcal{TH}$ has the homotopy-type of $S^1 \times S^1$.  
If $f \in \hat \K_{3,1}$ is a hyperbolic knot, the $O_2$-action preserves $\hat \K_{3,1}(f)$ if and only if the knot is invertible.  If the knot is invertible, there is an $O_2$-equivariant homotopy-equivalence $S^1 \times S^1 \to \hat \K_{3,1}(f)$.  The action of $O_2$ on $S^1 \times S^1$ is given by $A.(z_1,z_2) = (Az_1,Az_2)$ where $A \in O_2$.  Here $z_i \in S^1$ and $Az_i$ is the standard linear action of $O_2$ on $S^1$.  If the knot is not invertible, the component of the knot $f$ and its inverse $\overline{f}$ has the homotopy-type of $S^1 \times S^1 \times S^0$ and there is an $O_2$-equivariant homotopy-equivalence $S^1 \times S^1 \times S^0 \to \K_{3,1}(f) \cup \K_{3,1}(\overline{f})$ where the action of $O_2$ on $S^1 \times S^1 \times S^0$ is given by $A.(z_1,z_2,\epsilon) = (Az_1, Az_2, Det(A) \epsilon)$, where $\epsilon \in S^0 = \{\pm 1\}$. 
\end{itemize}
\begin{proof}
Both Brendle-Hatcher \cite{BH} and Theorem \ref{afflin}, have a central shrinking and linearization argument that assert that certain spaces of unlinks have the homotopy-type of the subspace consisting of linear embeddings.  Both arguments are highly analogous.  Although the Brendle-Hatcher argument is about the space of $k$-component unlinks (denoted by them as $\mathcal{AL}_{0,k}$), it applies equally well to the space of $(k+1)$-component KGLs (see Definition \ref{splicedef}), since the space of $(k+1)$-component KGLs fibre over the space of $k$-component unlinks.  One of the key theorems of Brendle-Hatcher is that $\mathcal{AL}_{0,k}$ has the homotopy-type of $\mathcal{R}_k$ (the subspace where all the circles are round i.e. geometric circles), moreover this space has a homotopy-equivalent subspace $\mathcal{SR}_k$ where the circles are `separated'.  The point being that elements of $\mathcal{SR}_k$ have a well-defined semi-linear ordering. One component $w_i$ is `lower' than another $w_j$ if the shell $S(w_i)$ bounds a ball containing $S(w_j)$.  Since the space of $(k+1)$-component KGLs fibres over the space of $k$-component unlinks, it is therefore homotopy-equivalent to the subspace where the underlying $k$-component unlink is separated.  The proof of proposition \ref{afflin} similarly gives a homotopy-equivalence between $\SP_{3,1}(k)$ and the subspace where $L_1, \cdots, L_k$ are separated.  Thus $\SP_{3,1}(k)$ has the homotopy-type of the space of $(k+1)$-component KGLs. 

Given $f \in \hat \K_{3,1}$, let $\hat\K_{3,1}(f)$ denote the path-component of $\hat\K_{3,1}$ containing $f$.  Let $C_f$ be the complement of an open tubular neighbourhood of the associated closed knot in $S^3$.  Let $Diff(C_f)$ denote the group of diffeomorphisms of $C_f$ which restrict to the identity on the boundary, then $\hat\K_{3,1}(f) \simeq BDiff(C_f)$.  This is a fairly standard argument based on the fact that the group of diffeomorphisms of the $3$-ball that fix the boundary point-wise, $Diff(D^3)$, is contractible \cite{Hatcher1, Hatcher2} (see \cite{Bud} or \cite{KS} for details on the homotopy-equivalence).  Let $T \subset C_f$ be the tori of the JSJ-decomposition of $C_f$. One can think of $T$ as defining a rooted tree (the `JSJ-tree' \cite{BudJSJ}) where the vertices are the path-components of $C_f$ split along $T$, and the edges are the path-components of $T$.  The root of the tree is the component of $C_f$ split along $T$ containing $\partial C_f$. Let $V$ consist of $C_f$ with the submanifold of $C_f$ corresponding to the leaves of the JSJ-tree removed. The complement of $V$ in $C_f$ is the union of disjoint non-trivial knot complements $\sqcup_{i=1}^k C_{f_i}$, where $f_i \in \hat \K_{3,1}$.  An observation that goes back to Schubert \cite{Sch2} (reproven in \cite{BudJSJ}) is that disjoint non-trivial knot complements in $S^3$ can be separated by disjoint embedded $3$-balls in $S^3$.  The operation of `unknotting' $f_1$ through $f_k$ gives a new embedding of $V$ in $S^3$ as the complement of an $(k+1)$-component link $\hat L \subset S^3$ corresponding to some $L \in \SP_{3,1}(k)$. The construction of $L$ can be made into a unique decomposition for $f$ provided we assert that $f$ is obtained by splicing i.e. $f$ is isotopic to $L.(f_1, \cdots, f_k)$.   Let $Diff(C_f,V)$ denote the subgroup of $Diff(C_f)$ which preserves $V$.  The inclusion $Diff(C_f,V) \to Diff(C_f)$ is known to be a homotopy-equivalence \cite{Hatcher1, Hatcher2} (see \cite{Bud, KS} for details). So we have a locally-trivial fibre bundle of topological groups $Diff(C_f,V) \to Diff(V)$. We use `locally trivial' in the sense common in the study of embedding spaces, that fibres can vary as one moves from component to component in the base, in particular they can be empty. The non-empty fibres can be identified with $\prod_{i=1}^k Diff(C_{f_i})$.  

Let $\SP_{3,1}(L)$ denote the path-component of $\SP_{3,1}$ corresponding to $L$. The re-embedding diffeomorphism $V \to C_L$ allows us to identify $Diff(V)$ with $Diff(C_L)$. Let $A_f$ be the maximal subgroup of $\Sigma_k \wr O_2$ preserving the component $\SP_{3,1}(L) \times \prod_{i=1}^k \hat\K_{3,1}(f_i)$ for the action of $\Sigma_k \wr O_2$ on $\SP_{3,1}(k)\times (\hat \K_{3,1})^k$ (see Proposition \ref{equivarianceprop}).  Applying the classifying-space functor to the locally-trivial fibre bundle of groups $Diff(C_f,V) \to Diff(C_L)$ gives a locally trivial fibre bundle with connected base space
$$ \prod_{i=1}^k \hat \K_{3,1}(f_i) \to \hat \K_{3,1}(f) \to \SP_{3,1}(L)/A_f.$$
By design the knots $f_i \in \mathcal{TH}$ for all $i$ (see Definition \ref{complexitydef}).  
The action of $\SP_{3,1}$ on $\hat\K_{3,1}$ gives us the central vertical map in a commuting diagram of onto fibrations
$$\xymatrix{  
 \prod_{i=1}^k \hat\K_{3,1}(f_i) \ar[r] \ar[d] & \SP_{3,1}(L) \times_{A_f} \prod_{i=1}^k \hat\K_{3,1}(f_i) \ar[r] \ar[d] & \SP_{3,1}(L)/{A_f} \ar[d] \\
 \prod_{i=1}^k \hat\K_{3,1}(f_i) \ar[r]        & \hat \K_{3,1}(f) \ar[r] &  \SP_{3,1}(L) / {A_f} }$$
Since the left and rightmost vertical arrows are homotopy-equivalences, the central vertical arrow is as well.  
For the claims describing the $O_2$-action on $\mathcal{TH}$, the key argument is to find suitable maximal-symmetry positions
for the closed versions of the knot in $S^3$. The equivariant maps to $\hat \K_{3,1}$ are given by a stereographic projection
construction which appears in detail in the proof of Theorem \ref{operadstruc}. 

Consider whether or not the splicing construction is an $O_2$-equivariant homotopy-equivalence.  By the $G$-Whitehead Theorem \cite{GWhit}, it would suffice to show that the map is a weak equivalence of $O_2$-spaces, meaning for every closed subgroup $H \subset O_2$, the splicing map is a homotopy-equivalence when restricted to the subspace fixed by $H$.
If a group $G$ acts on a space $X$ we denote the $G$-fixed point subspace of $X$ by $X^G$.   For $H$ any non-trivial closed subgroup of $SO_2$ this is immediate as only the linearly-embedded unknot is fixed by a non-trivial element of $SO_2$.  The only interesting case remaining is $H \simeq \Zed_2$, a subgroup whose fixed points $\K_{3,1}^H$ are knots in strong inversion positions. Stated another way, showing the splicing map from Theorem \ref{freenessthm} is an $O_2$-equivariant homotopy-equivalence amounts to showing that for strongly-invertible knots $f$, the space of strong inversion positions of $f$, $\K_{3,1}(f)^H$ is homotopy-equivalent to $\left(\SP_{3,1}(L) \times_{A_f} \prod_{i=1}^n \hat\K_{3,1}(f_i)\right)^H$, and the splicing map is such a homotopy-equivalence.  

Since 3-manifolds have equivariant JSJ-decompositions \cite{MS} and an equivariant Loop Theorem
\cite{JR}, the proof of Proposition 2.1 from \cite{BudJSJ} extends, giving the result that if a knot is isotopic to a non-trivial splice, and if that knot is strongly invertible, then one can put the knot into a position where it is simultaneouly strongly invertible and in the image of the splicing map.  Thus splicing gives an onto map 
$$\pi_0 \left(\SP_{3,1}(L) \times_{A_f} \prod_{i=1}^n \hat\K_{3,1}(f_i)\right)^H \to \pi_0 \left(\K_{3,1}(f)^H\right).$$
By the $\Zed_2$-equivariant isotopy extension theorem \cite{MK}, there is a fibre bundle $(Diff(C_f))^H \to (Diff(D^3))^H \to (\K_{3,1}(f))^H$.  By a standard cut-and-paste argument using \cite{Hatcher1} and \cite{Hatcher2}, $(Diff(D^3))^H$ is contractible.  Thus the space of strong invertibility positions for a knot has the homotopy-type of $B(Diff(C_f)^H)$. Repeating the above argument that $Diff(C_f)$ is a bundle over $Diff(C_L)$ with fibre $\prod_{i=1}^n Diff(C_{f_i})$ in this context, gives the homotopy-equivalence 
$$\left(\SP_{3,1}(L) \times_{A_f} \prod_{i=1}^n \hat\K_{3,1}(f_i)\right)^H \to \K_{3,1}(f)^H.$$
\end{proof}
\end{thm}

To make the statement of Theorem \ref{operadstruc} more compact, we introduce some terminology.  

\begin{defn}
Denote the $\Sigma^* \wr O_2$-suboperad of $\SP_{3,1}$ generated by the inclusion $\Cu_1' \subset \SP_{3,1}$ from Proposition \ref{cubesinclusion} be denoted $\overline{\Cu_1'}$.  Given a link $\hat L$ in $S^3$ call it {\it totally prime} if the JSJ-decomposition of $S^3 \setminus \hat L$ contains no manifolds diffeomorphic to the product of a circle and a punctured disc. Let $\mathcal{TP} \subset \SP_{3,1}$ be the subspace of $L \in \SP_{3,1}$ such that $\hat L$ is totally prime. 
\end{defn}

\begin{prop}\label{overlinec1prime}
As a $\Sigma^* \wr O_2$-operad $\overline{\Cu_1'}$ is equivalent to $\Cu_1' \rtimes O_2$.  The action of $O_2$ on $\Cu_1'$ factors through the homomorphism $O_2 \to \Zed_2$, $\Zed_2$ acting by mirror reflection
on $[-1,1]$. 
\begin{proof}
By design, $\overline{\Cu_1'}(k) = \Cu_1'(k) \times O_2^k$, as the act of taking the $\Sigma^* \wr O_2$-operadic closure of $\Cu_1'(k)$ amounts to adding all linear reparametrizations of the pucks. 
The result follows. 
\end{proof}
\end{prop}

Before proceeding to Theorem \ref{operadstruc}, we record some useful facts about cyclic and
dihedral groups acting on $S^3$.  For the next definition we will think of $\Zed_n \subset S^1 \subset \mathbb{C}$
as being the $n$-th roots of unity.  Given $p, q \in \Zed$ with $GCD(p,q)=1$, the $(p,q)$-embedding of
$\Zed_n$ in $SO_4$ is given by the action $\Zed_n \times \mathbb{C}^2 \to \mathbb{C}^2$ where $(z,(z_1,z_2)) \longmapsto (z^pz_1,z^qz_2)$. The {\it standard involution} of $S^3$ is the map $(z_1,z_2) \longmapsto (\overline{z_1},\overline{z_2})$. 

\begin{lem}\label{cyclicaction} 
Let $G$ be a finite subgroup of the group of orientation-preserving diffeomorphisms of $S^3$.  
Then $G$ is conjugate to a subgroup of $SO_4 \subset Diff^+(S^3)$.  If $G \subset SO_4$ is 
cyclic then it is conjugate to a $(p,q)$-action for some $p,q \in \Zed$ with $GCD(p,q)=1$. 
There is only one extension (up to conjugacy) of a $(p,q)$-action of $\Zed_n$ on $S^3$ to an
an action of $D_n$ on $S^3$.  If $n>2$ one of the involutions can be taken to be the standard involution.
When $n=2$ the extension of the $(0,1)$-action is by the antipodal map, as $D_2$ is abelian.
\begin{proof}
The fact that $G$ is conjugate to a subgroup of $SO_4$ is the `linearization' part of the elliptization conjecture
i.e. elliptization modulo the Poincar\'e conjecture. If $G$ acts freely, see \cite{MT}. If the action is not free,
see \cite{SmithConj}.  The remainder of this lemma can be derived by considering the eigenspaces of elements of $G$. 
\end{proof}
\end{lem}

Notice that the part of $S^3$ on which $G$ does not act freely has a rather simple structure. In the case that $G$ is cyclic it acts freely on $S^3$ if and only if $GCD(p,n) = GCD(q,n) = 1$.  If $GCD(p,n)=1$ but $GCD(q,n)>1$ there is the singular set $(\{0\}\times \mathbb{C}) \cap S^3$, which is a trivial knot.  If both $GCD(p,n)$ and $GCD(q,n)>1$ then the singular set is $\left((\mathbb{C} \times \{0\}) \cup (\{0\}\times \mathbb{C})\right) \cap S^3$, a Hopf link. In the case that $G$ is dihedral there are also the circles fixed by the involutions.

\begin{prop}\label{maxsym}
Let $L = (L_0,L_1,\cdots,L_k)$ be a hyperbolic link in $S^3$.  Then it has a maximal symmetry position with respect to the action of $\pi_0 Diff(S^3,L)$, meaning one can isotope $L$ into a position where the maps
$$Isom(S^3,L) \to \pi_0 Diff(S^3,L) \to Isom_{\mathcal{H}^3}(S^3\setminus L)$$
are isomorphisms. 
\begin{itemize}
\item $Isom(S^3,L)$ is the group of isometries of $S^3$ that preserve $L$ as a set -- there may or may not be fixed points on $L$. 
\item $\pi_0 Diff(S^3, L)$ is the link symmetry group i.e. the mapping class group of the pair $(S^3,L)$. 
\item $Isom_{\mathcal{H}^3}(S^3 \setminus L)$ is the group of hyperbolic isometries of the complement of $L$ which preserve meridional homology classes -- i.e. isometries of the link complement which admit continuous 
extensions $S^3\to S^3$.  
\item The map $\pi_0 Diff(S^3,L) \to Isom_{\mathcal{H}^3}(S^3\setminus L)$ is induced by Mostow rigidity -- i.e. restrict a diffeomorphism of the pair $(S^3,L)$ to $S^3 \setminus L$ and choose the unique hyperbolic isometry in that map's homotopy-class. 
\end{itemize}
If we demand that $(L_1,\cdots,L_k)$ is the trivial link, let $B_L$ denote the subgroup of 
$\pi_0 Diff(S^3,L)$ that preserves $L_0$ and the orientation of $S^3$.  Then one can isotope $L$ to ensure $B_L$ acts on $(S^3,L)$ by isometries of $S^3$, and where $L_i$ are {\it round circles} for all $i \in \{1,2,\cdots,k\}$, that is, the intersection of an affine 2-dimensional subspace of $\Real^4$ with $S^3$.
 
\begin{eg} For the Borromean rings $L$ on the left, $\pi_0 Diff(S^3,L)$ 
is the full octahedral group, having order $48$. An approximate maximal symmetry position on the left, and the maximal symmetry position for $B_L \simeq D_4$ and $L_1$, with
$L_2$ round is on the right.  In this picture the dotted blue circle is the singular set of the action of $\Zed_4$ on $S^3$. In the language of Lemma \ref{cyclicaction}, this is the $(2,1)$-action of $D_4$ on $S^3$ and the dotted blue
circle is $(\mathbb{C}\times\{0\}) \cap S^3$.  The dotted blue circle intersects both $L_1$ and $L_2$ in two points
each, but does not intersect $L_0$. 
{
\psfrag{l0}[tl][tl][1][0]{$L_0$}
\psfrag{l1}[tl][tl][0.8][0]{$L_1$}
\psfrag{l2}[tl][tl][0.8][0]{$L_2$}
$$\includegraphics[width=4cm]{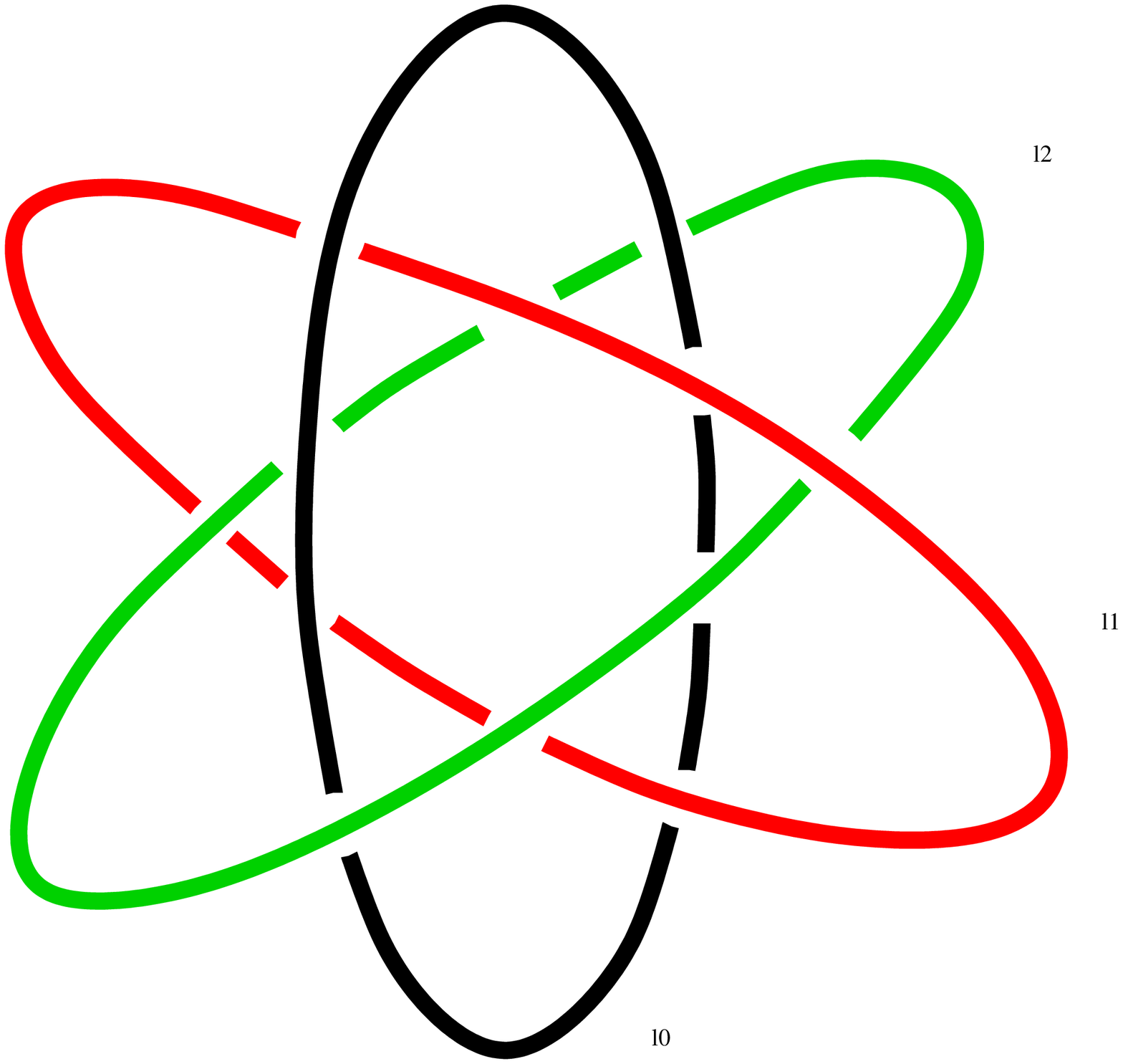} \hskip 2cm \includegraphics[width=4cm]{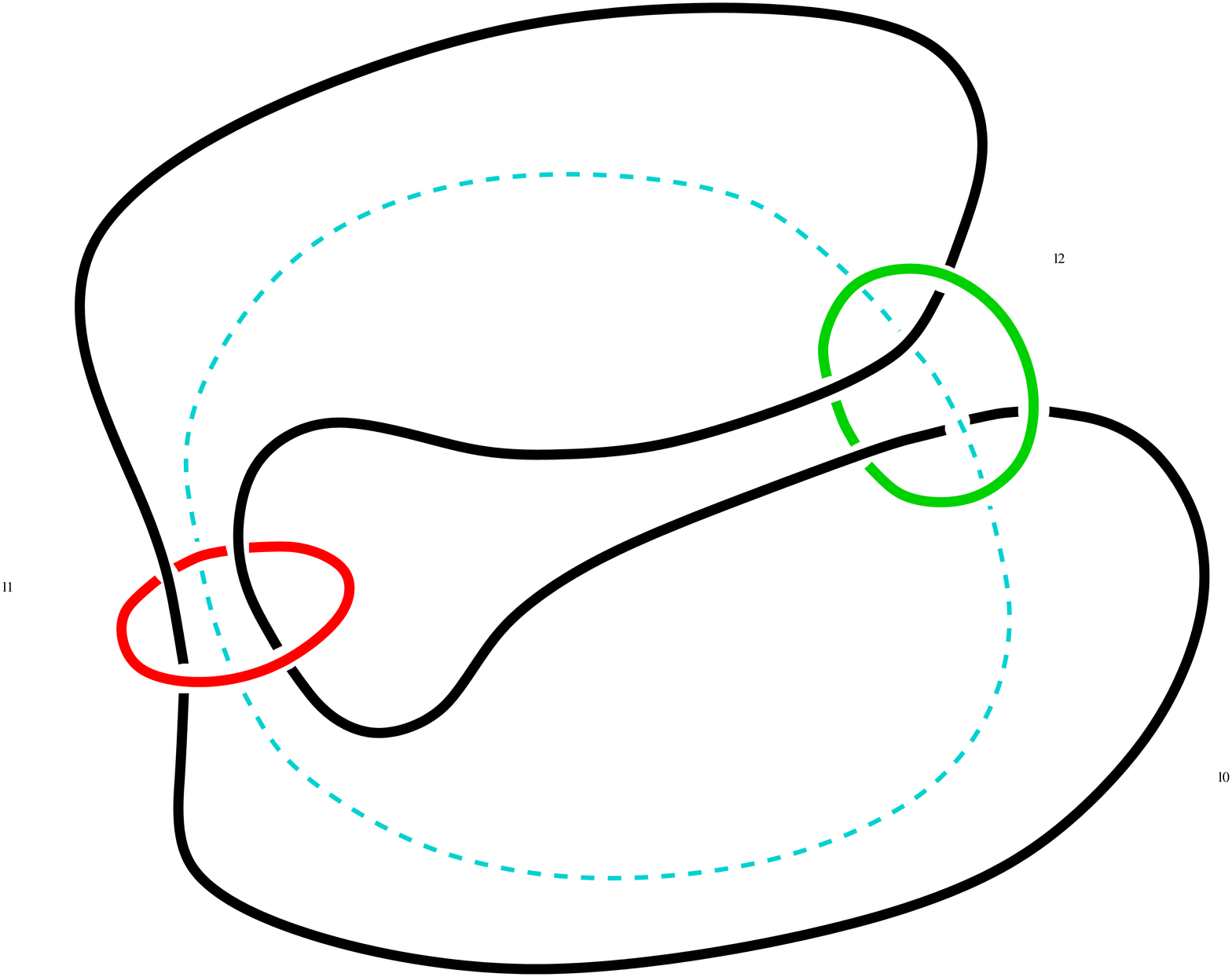}$$
}
\end{eg}
\begin{proof}
The existence of maximal symmetry positions is a standard amalgamation of several 
major theorems: 

\begin{itemize}
\item The group $Isom_{\mathcal{H}^3}(S^3 \setminus L)$ is finite, since isometry groups of complete finite volume hyperbolic 3-manifolds are finite. By definition, $Isom_{\mathcal{H}^3}(S^3 \setminus L)$ preserves the
longitudinal homology classes of $L$ so the action extends to an action of $Isom_{\mathcal{H}^3}(S^3 \setminus L)$ on $S^3$ giving an injective homomorphism $Isom_{\mathcal{H}^3}(S^3 \setminus L) \to Diff(S^3,L)$.
\item Due to the Elliptisation Theorem \cite{MT, SmithConj}, the action of 
$Isom_{\mathcal{H}^3}(S^3 \setminus L)$ on $S^3$ is conjugate to a linear action, i.e. there exists a diffeomorphism of $S^3$, $h : S^3 \to S^3$ such that the diagram commutes
$$\xymatrix{ Isom_{\mathcal{H}^3}(S^3 \setminus L) \times S^3 \ar[d]_{I \times h} \ar[r] & S^3 \\
             Isom_{\mathcal{H}^3}(S^3 \setminus h(L)) \times S^3 \ar[ur] & } $$
where the top horizontal arrow is the action of $Isom_{\mathcal{H}^3}(S^3 \setminus L)$ on $S^3$ and the diagonal arrow is a linear action of $Isom_{\mathcal{H}^3}(S^3 \setminus h(L))$ on $S^3$. $h(L)$ is the `maximal symmetry position' for $L$.  It is isotopic to $L$ since we can assume $h$ is orientation preserving. Orientation-preserving diffeomorphisms of $S^3$ are isotopic to the identity \cite{Cerf}.  
\item To complete the claim one uses work of Hatcher and Waldhausen that implies $Diff(S^3 \setminus L) \to HomEq(S^3 \setminus L)$ is a homotopy-equivalence, and by Mostow Rigidity that $Isom_{\mathcal{H}^3}(S^3 \setminus L) \to HomEq(S^3 \setminus L)$ is homotopy-equivalence, see Proposition 3.2 from \cite{KS} for details.  For the remainder
of the proof we replace $L$ with $h(L)$. 
\end{itemize}

To construct the maximal symmetry position for $B_L$, apply the Equivariant Sphere Theorem \cite{JR} of Jaco and Rubinstein to the $B_L$-manifold $S^3 \setminus \nu(L_1 \cup \cdots \cup L_k)$, where $\nu(L_1 \cup \cdots \cup L_k)$ indicates an open tubular neighbourhood of $L_1 \cup \cdots \cup L_k$ in $S^3$.  This gives us an equivariant collection $\mathcal S$ of embedded $S^2$'s in $S^3 \setminus \nu(L_1 \cup \cdots \cup L_k)$ which separate the manifold into a collection of punctured spheres ($S^3$) and punctured unknot complements ($S^1\times D^2$).  Think of $B_L$ as being a group of automorphisms of a rooted tree, the tree's vertices being the path-components of $S^3 \setminus \mathcal S$, and edges the path-components of $\mathcal S$.   Since finite groups acting on trees either fix a vertex or the centre of an edge, by replacing a sphere from $\mathcal S$ with the boundary of its equivariant tubular neighbourhood in $S^3$, we can arrange for there to be a vertex fixed by the action of $B_L$, i.e. some component of $S^3 \setminus \mathcal S$ is preserved by $B_L$. By Lemma \ref{cyclicaction} we have models for the action of the various stabilizers in $B_L$ on the components of $S^3 \setminus \mathcal S$. The components of $S^3 \setminus \mathcal S$ are punctured spheres so the action is the restriction of some $(p,q)$-embedding of a dihedral group in $SO_4$, in particular the action is linear.  Consider a component $B$ of $S^3 \setminus \mathcal S$ corresponding to a leaf of the tree, this is a 3-ball containing a single component of $L_1 \cup \cdots \cup L_k$. The subgroup of $B_L$ preserving $B$, if not trivial has singular set either an unknotted arc in $B$ or two unknotted arcs meeting at a central vertex.  Thus if $L_i$ is in $B$, $L_i$ either Hopf links the singular set or meets the singular set in two points. Either way, via a shrinking construction we can equivariantly linearize $L_i$ in $B$ to a round circle. This allows us to equivariantly shrink $B$ to the point that it is a small round ball.  Inductively, we can work from the leaves to the root of the tree associated to $\mathcal S \subset S^3$ and assume all the spheres and link components $L_1, \cdots, L_k$ are round.    By equivariant isotopy extension \cite{MK} we can isotope $L$ into a position such that $L_1, \cdots, L_k$ are round circles.
\end{proof}
\end{prop}

Theorem \ref{operadstruc} describes the equivariant homotopy-type of the operad $\SP_{3,1}$.  A key step in the argument is the construction of finite-dimensional subspaces of $\SP_{3,1}$ where the equivariant homotopy type is explicitly understood.  The most elaborate case consists of the components of $\SP_{3,1}$ containing hyperbolic links.  Using the maximal symmetry positions of hyperbolic links, via a stereographic projection construction we will create these finite-dimensional families in $\SP_{3,1}$, allowing us to understand the $\Sigma^*_k \wr O_2$-equivariant homotopy-type of $\SP_{3,1}(k)$.  Will need some conventions relating the group $\Sigma^*_k \wr O_2$ to the geometry of the link $\hat L = (\hat L_0, \cdots, \hat L_k) \subset S^3$. 

\begin{defn}\label{symconventions}
Define 
$$F\hat L = F\hat L_0 \times \Sigma_k \times \prod_{i=1}^k UT\hat L_i$$
 where $UT\hat L_i$ is the unit tangent bundle to $\hat L_i$, and $FL_0$ is the frame bundle of $L_0$, meaning 
$$F\hat L_0 = \{ (p,w_1,w_2) : p \in L_0, w_1 \in T_p L_0, w_2 \in \Real^4 \text{ and the triple } (p, w_1, w_2) \text{ is orthonormal}\}.$$  
$F\hat L$ should be thought of as the minimal data to uniquely describe:
\begin{itemize}
\item a constant-speed diffeomorphism $S^1 \to \hat L_0$, 
\item a constant-speed diffeomorphism $\sqcup_k S^1 \to \hat L_1 \cup \cdots \cup \hat L_k$
\item a unit-length normal vector field to $\hat L_0$ for which its covariant derivative is parallel along $\hat L_0$, moreover we demand this normal vector field does not homologically link $\hat L_0$.  Here `parallel' means with respect to the connection on the normal bundle induced by orthogonal projection. 
\end{itemize}
By design there is a left action of $B_L$ on $F\hat L$ given by post-composition of these parametrizations with an isometry of $S^3$.  There is also a right action of $Aut(\nu S^1) \times \Sigma_k \wr O_2$ on $F \hat L$ given by pre-composition with an isometry of the parametrizing domain $\nu S^1 \sqcup (\sqcup_k S^1)$, moreover these two actions on $F\hat L$ commute.   We use the convention that $\nu S^1$ is the trivial $S^1$-bundle over $S^1$, and $Aut(\nu S^1) \equiv (S^1 \times S^1) \rtimes \Zed_2$ is automorphisms of the bundle that are orientation-preserving on the total space.   Since any two parametrizations differ by precomposition with an element of $Aut(\nu S^1) \times \Sigma_k \wr O_2$, $F\hat L$ is an $Aut(\nu S^1) \times \Sigma_k \wr O_2$-torsor.  This induces a canonical injection $B_L \to Aut(\nu S^1) \times \Sigma_k \wr O_2$.  The composition with the projection $B_L \to Aut(\nu S^1) \times \Sigma_k \wr O_2 \to Aut(\nu S^1)$ is an embedding of groups. 
\end{defn}



\begin{eg} A hyperbolic link $L$ with $B_L \simeq D_3$ (dihedral group of triangle).  
{
\psfrag{l0}[tl][tl][1][0]{$L_0$}
\psfrag{l1}[tl][tl][0.8][0]{$L_1$}
\psfrag{l2}[tl][tl][0.8][0]{$L_2$}
$$\includegraphics[width=4cm]{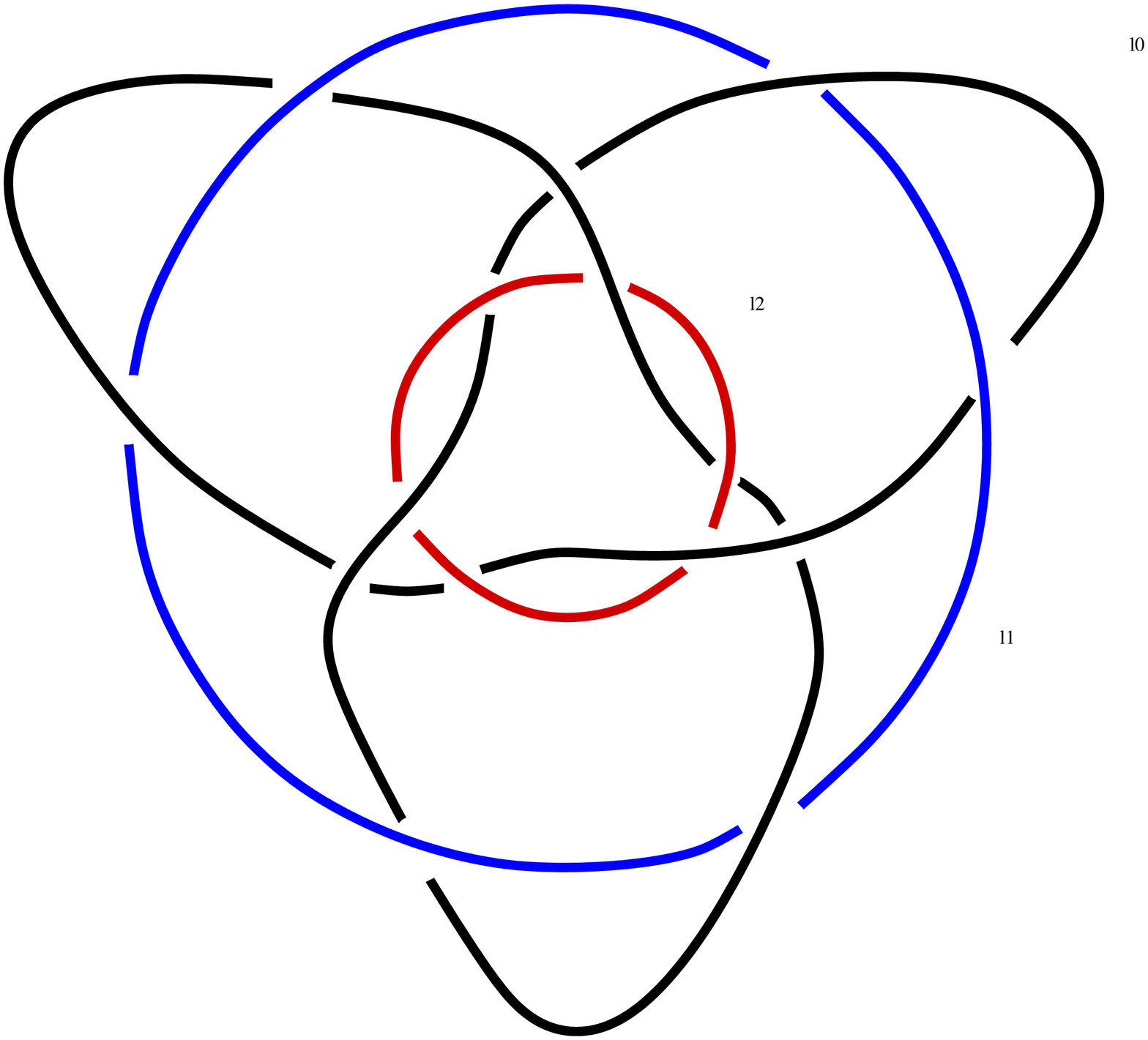}$$
}
\end{eg}



A key geometric construction in Theorem \ref{operadstruc} is a map from $F\hat L$ to the space of KGL's (see Definition \ref{splicedef}). Given a point of $W \in F\hat L$, it determines a point of $F\hat L_0$.  This is a point on $L_0$ together with a unit tangent vector and a unit normal vector.  We think of $S^3$ as the unit vectors in $\Real^4$.  Since all three vectors are orthogonal in $\Real^4$, if we think of them as column vectors of a $4 \times 3$-matrix, they extend uniquely to an element $A_W \in SO_4$.   Let $\sigma_W \in \Sigma_k$ be the permutation specified by $W \in F\hat L$. 

\begin{eg} The Whitehead link in its maximal symmetry position in $S^3$ with $\hat L_1$ a round circle, together with a sampling of stereographic projections along $\hat L_0$. 
{
\psfrag{L0}[tl][tl][1.1][0]{$\hat L_0$}
\psfrag{L1}[tl][tl][0.8][0]{$\hat L_1$}
\psfrag{l2}[tl][tl][0.8][0]{$L_2$}
\psfrag{l3}[tl][tl][0.8][0]{$L_3$}
$$\includegraphics[width=14cm]{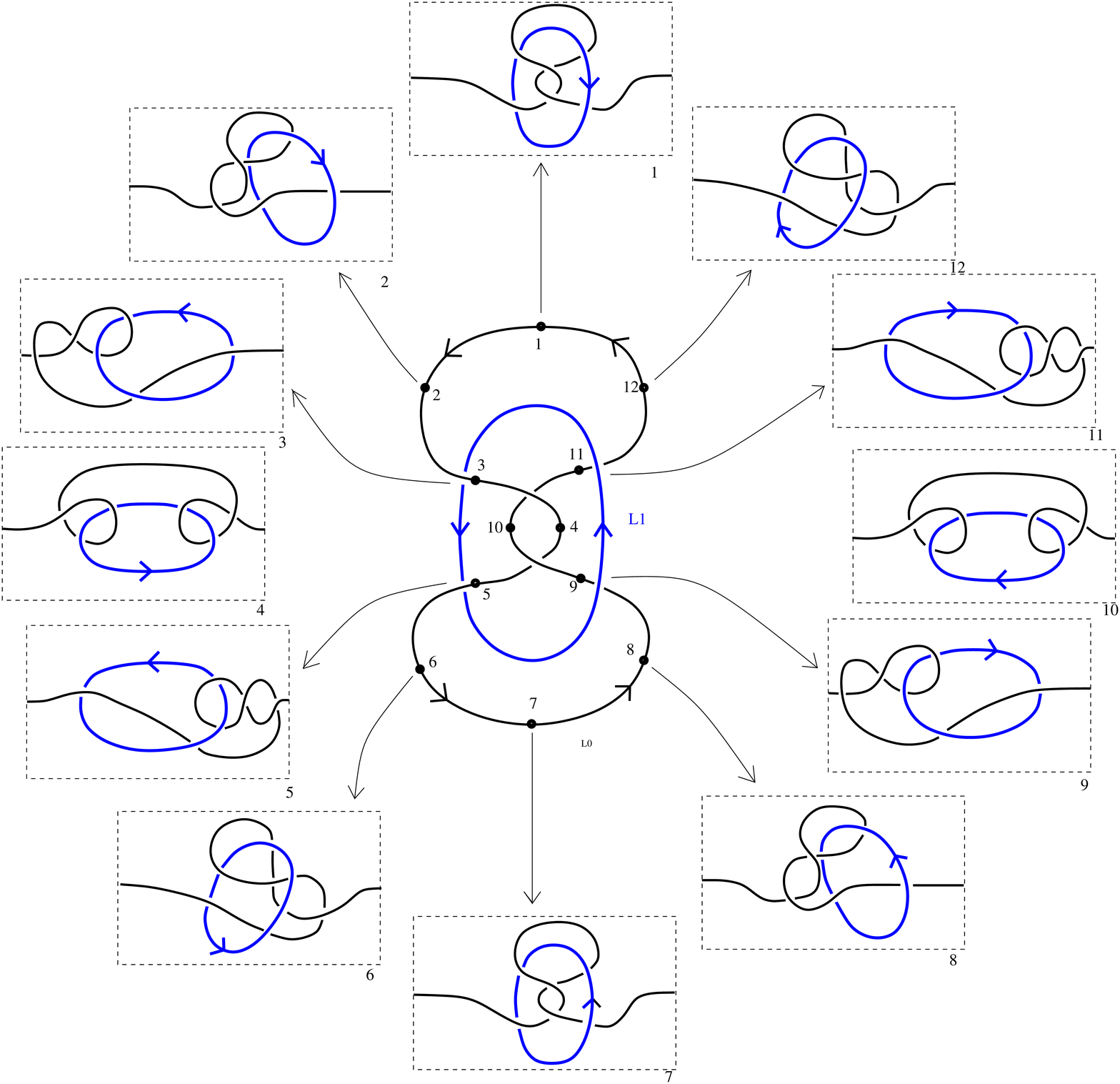}$$
}
\end{eg}

Let $f_{i,W} : S^1 \to \hat L_{\sigma(i)}$ be the uniquely-prescribed constant-speed parametrization specified by $W$, i.e. such that the derivative of $f_{i,W}$ at $1$ in the counter-clockwise direction is (up to a positive multiple) the unit tangent vector of $\hat L_{\sigma(i)}$ specified by $W$. Consider stereographic projection to be a map $p_a:S^n \to T_a S^n$ for any $a\in S^n$.  Conjugation of $A_W^{-1} f_{0,W}$ by stereographic projection
$$\xymatrix{S^1 \ar[r]^{A_W^{-1}f_{0,W}} \ar[d]^{p_1} & S^3 \ar[d]^{p_1} \\
            \Real \ar[r]      & \Real^3 }$$
produces a map $p_1 A_W^{-1}f_{0,W} p_1^{-1}$ which is `almost' an element of $\K_{3,1}$.  Similarly, composing $A_W^{-1} f_{i,W}$ with $p_1 : S^3 \to \Real^3$ produces the collection $(p_1 A_W^{-1}f_{0,W} p_1^{-1}, p_1 A_W^{-1} f_{1,W}, \-\cdots,\- p_1 A_W^{-1} f_{k,W})$ which is `almost' a KGL.  This collection is an embedding $\Real \cup (\sqcup_k S^1) \to \Real^3 \equiv T_1S^3$ which fails to be a KGL precisely when $f_{0,W}$ fails to be linear in a sufficiently large neighbourhood of $1$, or if $\hat L_1, \cdots, \hat L_k$ get too close to $\hat L_0$ in the sense that their stereographic projections may not be contained in $I \times D^2$ (see Definition \ref{splicedef}).  This is not a serious obstacle in that we can equivariantly linearize $f_{0,W}$ near $1$ and suitably rescale via a hyperbolic transformation of $S^3$ at $\pi(v)$, at which point stereographic projection will give an actual KGL.  A key point in this argument is that since stereographic projection preserves round circles, the stereographic projections of $\hat L_1, \cdots, \hat L_k$ are round circles in $\Real^3$, so they bound canonical flat discs which can be fattened into hockey pucks. 

\begin{thm}\label{operadstruc}
As an $\Sigma^*\wr O_2$-operad, $\SP_{3,1}$ is the free product of $\overline{\Cu_1'}$ and $\mathcal{TP}$. Moreover, $\mathcal{TP}$ is a free $\Sigma^* \wr O_2$-operad, freely generated by the subspace $\left(\sqcup_{k \in \Nat} \mathcal{HGL}_k\right) \sqcup \mathcal{SFL}$ of $\SP_{3,1}$ consisting of:
\begin{enumerate}
\item $\mathcal{SFL}$, these are the $2$-component Seifert links (Hopf link not included) in $\SP_{3,1}(1)$. Once we close these links to be links in $S^3$, these are the links $\mathcal{S}^{(p,q)}$ from \cite{BudJSJ} with $(p,q)\in \Zed^2$, $GCD(p,q)=1$ and $p \nmid q$, i.e. the Seifert link of type $(p,q)$ is a $2$-component link in $S^3$ consisting of two fibres in a $(p,q)$-Seifert fibring of $S^3$, one fibre singular, the other not.  $\mathcal{S}^{(p,q)} = (\{ (z_1,z_2) \in \mathbb C^2 : z_1^p=z_2^q\} \cap S^3) \cup (S^1 \times \{0\}) \subset S^3$. 
\item Hyperbolic links $k \in \{1,2,3,\cdots\}$, meaning that $L \in \SP_{3,1}(k)$ belongs to $\mathcal{HGL}_k$ if and only if the complement of the corresponding closed link $\hat L$ in $S^3$ has a complete hyperbolic structure of finite-volume. 
\end{enumerate}
Restating the above in a different formalism, if we restrict the structure map 
$$\SP_{3,1}(k) \times_{\Sigma_k \wr O_2} \prod_{i=1}^k \SP_{3,1}(j_i) \to \SP_{3,1}(\sum_{i=1}^k j_i)$$
to the appropriate path-components of the domain and range respectively corresponding to a non-redundant splice, then it is an $( \Sigma^*_{(j_1+\cdots+j_k)} \wr O_2 )$-equivariant homotopy-equivalence between those components.
\begin{itemize}
\item[(1)] $\mathcal{SFL}$ has the homotopy-type of a disjoint union of countably-many tori $S^1 \times S^1$, two for every Seifert link $\mathcal{S}^{(p,q)}$.  Let $\SP_{3,1}(\mathcal{S}^{(p,q)})$ denote the path-component of $\SP_{3,1}$ corresponding to the $(p,q)$-Seifert link $\mathcal{S}^{(p,q)}$. Only an index two subgroup of $\Sigma_1^* \wr O_2 = O_2 \times O_2$ preserves $\SP_{3,1}(\mathcal{S}^{(p,q)})$, so let $(\Sigma_1^* \wr O_2).\SP_{3,1}(\mathcal{S}^{(p,q)})$ be the union of the path-components of $\SP_{3,1}$ containing the link corresponding to $\mathcal{S}^{(p,q)}$ and its image under the action of $\Sigma_1^* \wr O_2$. There is $O_2^2$-equivariant homotopy-equivalence $S^1 \times S^1 \times S^0 \to (\Sigma_1^* \wr O_2).\SP_{3,1}(\mathcal{S}^{(p,q)})$, where the right-action of $O_2^2$ on $S^1 \times S^1 \times S^0$ given by $(z_1,z_2,\epsilon).(A_1,A_2) = (A_1^{-1}z_1,A_2^{-1}z_2, \epsilon Det(A_1A_2))$. 
\item[(2)] $\mathcal{HGL}_k$ has the homotopy-type of a disjoint union of a countable collection of tori of the form $(S^1 \times S^1) \times (S^1)^k$.  Notice that the action of $\Sigma^*_k \wr O_2$ may permute path-components of $\mathcal{HGL}_k$, depending on the symmetry properties of $L \in \mathcal{HGL}_k$.   So we will describe the homotopy type of not just one path component, but the union of all path components containing $L$ and its image under the action of $\Sigma^*_k \wr O_2$. Denote this subspace of $\mathcal{HGL}_k$ by $(\Sigma^*_k \wr O_2).\mathcal{HGL}_k(L)$. Let $\hat L$ be the associated closed link in $S^3$, in its maximal symmetry position (Proposition \ref{maxsym}).  There is a $\Sigma^*_k \wr O_2$-equivariant homotopy-equivalence
$$ \Pi : F \hat L / B_L \to (\Sigma^*_k \wr O_2).\mathcal{HGL}_k(L).$$
where the action of $\Sigma^*_k \wr O_2$ on $F \hat L$ is given by considering $F\hat L$ as an $Aut(\nu S^1) \times \Sigma_k \wr O_2$-torsor, see Definition \ref{symconventions}. $(\{1\} \times S^1) \rtimes \Zed_2 \equiv O_2$ is a subgroup of $Aut(\nu S^1) = (S^1 \times S^1) \rtimes \Zed_2$. Think of $S^1 \times \{1\} \times \{0\}$ as the subgroup of $Aut(\nu S^1)$ corresponding to pure translational reparametrizations of $\hat L_0$, so the inclusion $O_2 \to Aut(\nu S^1)$ corresponds to reversing the tangent vector to $\hat L_0$ and rotating the normal vector.  Our inclusion $\Sigma_k^* \wr O_2 = O_2 \times \Sigma_k \wr O_2 \to Aut(\nu S^1) \times \Sigma_k \wr O_2$ is induced by this inclusion. 
\end{itemize}
\begin{proof}
The up-to-isotopy uniqueness statement for the splice decomposition was given in \cite{BudJSJ}.  That the splicing map restricts to an equivariant homotopy-equivalence for non-redundant splices, this argument is essentially the same as the proof of Theorem \ref{freenessthm}, with little modification beyond what is explained below. 

The homotopy-types of the spaces $\mathcal{SFL}$ and $\mathcal{HGL}_k$ are described in \cite{KS}, although the maps provided in that paper do not respect the $\Sigma^* \wr O_2$-action. Below we give a short summary of how the $\Sigma^* \wr O_2$-equivariant homotopy-type of each component of $\SP_{3,1}$ are determined.

The component of $\SP_{3,1}(k)/\Sigma_k \wr O_2$ corresponding to $L$ has the homotopy-type of the classifying space of a group of diffeomorphisms of a manifold $C_L$, denoted $Diff(C_L)$.  $C_L$ is the complement of an open tubular neighbourhood of $\hat L$ in $S^3$. $Diff(C_L)$ denotes the group of diffeomorphisms of $C_L$ which restrict to the identity on the boundary-component of $C_L$ corresponding to $\hat L_0$.  We also require the diffeomorphisms to preserve the homology classes (up to sign) of the set of meridians corresponding to $\hat L_1, \cdots, \hat L_k$ respectively, as this ensures the diffeomorphisms of $C_L$ extend to diffeomorphisms of $S^3$.  

In the case $C_L$ is Seifert-fibred, the diffeomorphism group has the homotopy-type of the fibre-preserving subgroup \cite{Hatcher1}. 
\begin{itemize}
\item For a Keychain link this group has the homotopy-type of the braid group on $k$ strands. 
\item For a Seifert link, $k=1$ and it has the homotopy-type of $\Zed$. The generator is a meridional Dehn twist about a torus of $\partial C_L$ corresponding to $\hat L_0$.  
\end{itemize}
In the hyperbolic case, Proposition \ref{maxsym} demonstrates that the full group of diffeomorphisms of $C_L$ has the homotopy-type of the group of hyperbolic isometries of $C_L \setminus \partial C_L \equiv S^3 \setminus \hat L$. The subgroup of $Isom_{\mathcal{H}^3}(S^3\setminus \hat L)$ that preserves the $\hat L_0$ cusp acts faithfully that cusp, so the restriction map from the diffeomorphism group of $C_L$ that preserves the boundary torus corresponding to $\hat L_0$ to the diffeomorphism group of that torus gives us an extension
$$0 \to \Zed^{2} \to \pi_1 (\SP_{3,1}(L) / B_L \wr O_2) \to F \to 0$$
where $F$ is a finite cyclic group with at most one generator by the `No Bad Monodromy' result \cite{KS}. $F$ can be understood as the translational symmetries of $\hat L_0$ induced by elements of $Isom_{\mathcal{H}^3}(S^3\setminus \hat L)$ from the perspective of Definition \ref{symconventions}. The $\Zed^2$ kernel consists of all the Dehn twists about a torus in the interior of $C_L$ which are parallel to the boundary torus corresponding to $\hat L_0$.  The extension is non-split provided $F$ is non-trivial. This is because the solution to the extension problem are `fractional Dehn twists' \cite{KS}.  This means that a diffeomorphism of $C_L$ that fixes the boundary torus (corresponding to $\hat L_0$) pointwise,  can be isotoped to agree with an isometry of $S^3 \setminus \hat L$ away from a collar neighbourhood of the fixed torus.  Inside that collar neighbourhood the diffeomorphism is free to be arbitrary translations of the torus fibers.   $\pi_1 (\SP_{3,1}(L) / B_L \wr O_2)$ is therefore free abelian of rank two. 

Since $\SP_{3,1}(k)$ fibres over the space of KGLs, the remainder of the proof is devoted to constructing an equivariant lift of the stereographic projection construction following Definition \ref{symconventions} to a $\Sigma_k^* \wr O_2$-equivariant map $F\hat L \to (\Sigma_k^* \wr O_2).\mathcal{HGL}_k(L)$ which descends to a homotopy-equivalence $F\hat L/B_L \to (\Sigma_k^* \wr O_2).\mathcal{HGL}_k(L)$.
 
To begin, we need to `fatten' $\hat L_0$, i.e. choose a $B_L$-equivariant tubular neighbourhood $\Upsilon$ of $\hat L_0$ in $S^3$ \cite{MK}.  Let $\nu_\epsilon S^1 = \{(z_1,z_2) \in \mathbb{C}^2 : |z_2| \leq \epsilon \} \cap S^3$ for any $0 < \epsilon < 1$, considering it to be the total-space of bundle over $S^1$ via the projection map $(z_1,z_2)\longmapsto z_1 \in S^1$. Trivialize the $B_L$-equivariant tubular neighbourhood explicitly, considering the trivialization to be a fibre-preserving diffeomorphism $\omega : \nu_\epsilon S^1 \to \Upsilon$.  If $\epsilon$ is sufficiently small, we can ensure $D\omega$ is conformal-linear along $S^1 \times \{0\}$, and by choosing the constant-speed parametrization of $\hat L_0$ we can ensure the conformal factor is constant.  The partial derivative of $\omega$ at $(z_1,z_2)$ in the direction of $(0,z_1)$ is a normal vector field along $\hat L_0$, and as in Definition \ref{symconventions}, we can choose it so that its covariant derivative is parallel along $\hat L_0$, and it does not homologically link with $\hat L_0$.  This reduces our choice of $\omega$ to the choice of $\epsilon$ and the initial data in $F\hat L_0$. 

Let $A : F \hat L \to F\hat L_0$, $\sigma : F \hat L \to \Sigma_k$, $u_i : F \hat L \to UT\hat L_i$ and $p_i : F\hat L \to \hat L_i$ be projection maps for $i \in \{1,2,\cdots, k\}$, for the product $F \hat L = F \hat L_0 \times \Sigma_k \times \prod_{i=1}^k UT\hat L_i$.  Using the conventions of Definition \ref{symconventions}, we will choose to think of $A$ as a map $A : F\hat L \to F\hat L_0 \equiv SO_4$. Given $W \in F\hat L$ let $\omega_W : \nu_\epsilon S^1 \to \Upsilon$ be the precomposition of $\omega$ with the uniquely-determined rigid motion $\nu_\epsilon S^1 \to \nu_\epsilon S^1$ so that $\omega_W$ agrees with some positive multiple of $A_W$ to first order at $(1,0) \in \nu_\epsilon S^1$. Let $g_W$ be the unique hyperbolic conformal transformation of $S^3$ fixing $\omega_W(1,0)$ such that $D(g_W \circ \omega_W)_{(1,0)} : T_{(1,0)}S^3 \to T_{\omega_W(1,0)}S^3$ is an isometry, therefore equal to $A_W$. $A_W^{-1} \circ g_W \circ \omega_W$ fixes $(1,0)$ and its derivative is the identity on $T_{(1,0)} S^3$.


Next we will apply a local linearization process to the embedding $A_W^{-1} \circ g_W \circ \omega_W$ at $(1,0)$.  Before that, a small digression into two standard linearization processes and how they can be related.  Given a diffeomorphism $f : U \to V$ where $U, V\subset \Real^n$ are open subsets of Euclidean space such that $0 \in U$, $f(0)=0$ and $Df_0$ is the identity $Df_0 = Id_{\Real^n}$, the {\it rescaling linearization process} means the homotopy $F_t(x) = \frac{1}{1-t}f((1-t)x)$.  Notice that at time $t$, the domain of $F_t$ is $\frac{1}{1-t}U$, and the image of $F_t$ is $\frac{1}{1-t}V$. We can extend $F$ to $t=1$ by $F_1(x) = x$ for all $x \in \Real^n$.  This is a variant of what is sometimes called the {\it Alexander Trick}.   Notice that if $Df$ has a Lipschitz constant $||Df_x - Df_y|| \leq K|x-y|$ for $x,y \in U$, then the Lipschitz constant for $D(F_t)$ is $(1-t)K$.   Similarly, the Lipschitz constant for the Hessian of $F_t$ is $(1-t)^2$ times the Lipschitz constant for the Hessian of $f$.  The second linearization process we consider is the {\it straight-line homotopy}.  More precisely, consider the problem of asking when the straight-line homotopy $G_t$ between $f$ and $Id_U$ is an isotopy: $G_t(x) = (1-t)f(x) + tx$.  One can check, in order for this to be an isotopy, it is sufficient for the Lipschitz constants for $Df$ and the Hessian of $f$ to be sufficiently small over the domain of $f$.  Combining the two linearization processes, we can say that given a ball neighbourhood $U'$ of $0 \in U$ such that $\overline{U'} \subset U$, there is some $\epsilon$, $0 < \epsilon < 1$ such that the straight-line homotopy between $F_\epsilon$ on $U'$ and the identity map $Id_{U'}$ is an isotopy.  Moreover, $\epsilon$ can be chosen to depend smoothly on the $C^2$-norm of $f$. 

Our linearization process for $A_W^{-1} \circ g_W \circ \omega_W$ is similar to what we did with $f : U \to V$. Multiplication by $(1-t)$ does not make sense on $S^3$, so we conjugate by the hyperbolic conformal transformations that fix the point $(1,0) \in S^3$.  These conformal transformation conjugate via stereographic projection $(1,0) \in S^3$ to multiplication by $(1-t)$ in $T_{(1,0)} S^3$ so they are completely analogous.  Specifically, given $p \in S^n$ and $t \in (0,\infty)$ the map $M_{p,t} : S^n \to S^n$ is multiplication by $t$ in $T_p S^n$ conjugated by stereographic projection at $p$ to be a map $S^n \to S^n$ fixing $p$ and $-p$. 

$$\xymatrix@R=7pt{ S^n \ar[r]^{M_{p,t}} \ar@{}[d]|{\upepsilon} & S^n \ar@{}[d]|{\upepsilon} \\
                    q  \ar@{|->}[r] & \frac{\left((t-1)^2(q \cdot p)+t^2-1\right)p + 2tq}{(t^2-1)(q \cdot p)+t^2+1} }$$

$q \cdot p$ denotes the standard Euclidean inner product $S^n \subset \Real^{n+1}$.   Our linearization process will start by considering the family of embeddings for $t \in [0,1)$

$$g_{W,t} = M_{(1,0), (1-t)}^{-1} \circ A_W^{-1} \circ g_W \circ \omega_W \circ M_{(1,0), (1-t)}$$

Notice when $t=0$, $g_{W,0} = A_W^{-1} \circ g_W \circ \omega_W$. The domain of $g_{W,t}$ is $M_{(1,0),(1-t)}^{-1} (\nu_\epsilon S^1)$, which for $t \in [0,1) \cap [1- \frac{\epsilon}{\sqrt{1-\epsilon^2}},1)$ contains the {\it right hemi-sphere} $HR = \{ (x,y,z,w) \in S^3 \subset \Real^4: x \geq 0 \}$. Moreover, for $t$ sufficiently close to $1$, $g_{W,t}$ approximates the identity map on the right hemi-sphere, {\it uniformly} in the $C^2$-topology.  Thus for $t$ sufficiently large we can ensure the geodesic/straight-line homotopy in $S^3$ from $g_{W,t | HR}$ to the identity map on the right hemi-sphere is an isotopy.  The equivariant isotopy extension theorem \cite{MK} allows us to extend this linearization of $g_{W,t|HR}$ to and isotopy of $g_{W,t}$.  Let $\Omega_W$ denote the resulting embedding which is linear on $HR$. It is almost never the case $\nu_\epsilon S^1$ is contained in the domain of $\Omega_W$, but using the convensions of Definition \ref{symconventions}, $p_1 \circ \Omega_W \circ p_1^{-1}$ is defined on a neighbourhood of $\Real \times \{0\}$ in $\Real^3$, moreover, it is the identity on a neighbourhood of $(\Real \setminus (-1,1)) \times D^2$. So the composite $p_1 \circ \Omega_W \circ p_1^{-1} \circ R_h$ is defined and is an element of $\hat \K_{3,1}$ for some $h \in (0,1]$ (see Proposition \ref{shrinkingprop} for the definition of $R_h$).  Moreover, we can choose $h$ to vary continuously with $W \in F\hat L$. Denote this element $f_{0,W} \in \hat K_{3,1}$. 

Given $i \in \{1,2,\cdots,k\}$ let $\overline{u_i(W)} : S^1 \to \hat L_i$ be the constant-speed parametrization of $\hat L_i$ such that $\overline{u_i(W)}(1) = p_i(W)$, and the derivative of $\overline{u_i(W)}$ at $1$ in the counter-clockwise direction is a positive multiple of $u_i(W) \in UT\hat L_i$. Notice that $M^{-1}_{(1,0),(1-t)} \circ A^{-1}_W \circ \overline{u_i(W)}$ is a parametrization of a round circle in $[-1,1] \times D^2$ and so it bounds a flat disc.  We choose a thickening of that disc and define it to be $f_{\sigma^{-1}_W(i),W} : [-1,1]\times D^2 \to [-1,1]\times D^2$.   Putting these all together, we have a map $F\hat L \to (\Sigma_k^* \wr O_2).\mathcal{HGL}_k(L)$ which is {\it by design} $\Sigma^*_k \wr O_2$-equivariant.  Since our action of $B_L$ on $F \hat L$ commutes with the action of $\Sigma^*_k \wr O_2$, this map descends to a $\Sigma^*_k \wr O_2$-equivariant map
$$F\hat L/B_L \to (\Sigma_k^* \wr O_2).\mathcal{HGL}_k(L)$$
for which we can check is an equivariant homotopy-equivalence.  Given $W \in F \hat L$, $A_W$ is a framed point in $\hat L_0$, denote this point by $p_0(W)$.  Notice that the action of $\Sigma^*_k \wr O_2$ leave the points $p_0(W)$ fixed.  Thus, the only fixed points of the action of $\Sigma_k^* \wr O_2$ on $F\hat L/B_L$ come from isometries of $S^3$ which {\it reverse} the orientation of $\hat L_0$.  So as in the proof of Theorem \ref{freenessthm} we can use the fact that both $\Sigma_k \wr O_2 \to F\hat L/B_L \to (F\hat L/B_L)/(\Sigma_k \wr O_2)$ and $\Sigma_k \wr O_2 \to (\Sigma_k^* \wr O_2).\mathcal{HGL}_k(L) \to \left((\Sigma_k^* \wr O_2).\mathcal{HGL}_k(L) \right)/(\Sigma_k \wr O_2)$ are principal $\Sigma_k \wr O_2$-bundles. The remaining $SO_2$-actions are free on both, there is only various order two mirror reflection subgroups of $O_2$ that have fixed-points, if the link is strongly invertible (for hyperbolic links, strong invertibility is implied by invertibility).  The homotopy-type of the space of strongly-invertible positions in $\left((\Sigma_k^* \wr O_2).\mathcal{HGL}_k(L) \right)/(\Sigma_k \wr O_2)$ is computed just as in Theorem \ref{freenessthm}, where we see immediately it has the homotopy-type of a product of two circles and our stereographic projection construction is by design an equivariant homotopy-equivalence. 
\end{proof}
\end{thm}

\begin{cor}
$\SP_{3,1}$ contains a homotopy-equivalent suboperad such that each component is finite-dimensional.
\end{cor}

The finite-dimensional suboperad of $\SP_{3,1}$ is of course the operad freely-generated by the suboperad $\overline{\Cu_1'} \equiv \Cu_1' \rtimes O_2 \subset \SP_{3,1}$ and the images of the maps $\Pi : F\hat L / B_L \to \SP_{3,1}$ for the hyperbolic links $L$ in the splicing operad. 

Theorems \ref{operadstruc} and \ref{freenessthm} say that the $O_2$-equivariant homotopy-type of $\hat K_{3,1}$, and the $\Sigma^* \wr O_2$-equivariant homotopy-type of $\mathcal{SP}_{3,1}$ is completely prescribed by the action of $B_L$ on $F \hat L$ for the hyperbolic links $L$ in the splicing operad. Determining which such representations arise is called the {\it realization problem} for the space of knots $\K_{3,1}$.  If one was only interested in the homotopy-type of $\K_{3,1}$ and $\mathcal{SP}_{3,1}$ respectively, one could take $\pi_0$ of these actions, and consider it to be a homomorphism $B_L \to \Sigma^*_k \wr \Zed_2$, and ask which such representations arise?  This is another variant of the realization problem. 

The next proposition points out that Proposition \ref{maxsym} gives new restrictions on which such representations can occur.  For the purpose of the realization problem a representation $\Zed \to \Sigma_k \wr \Zed_2$ is only interesting up to conjugacy.  Conjugacy classes in the symmetric group are traditionally specified by cycle decompositions, which are essentially partitions of the set $\{1,2,\cdots,k\}$.  The group $\Sigma_k \wr \Zed_2$ should be thought of as the signed permutation group, and conjugacy classes have a {\it signed cycle decomposition}. A signed cycle that preserves all signs is denoted $(a_1, a_2, \cdots, a_j)$. Let `$(a_1, a_2, \cdots, a_j)-$' denote the signed cycle type $a_1 \to a_2 \to \cdots \to a_j \to -a_1$, meaning all signs are preserved except the last one, which reverses sign. Thus $(a_1, a_2, \cdots, a_j)-$ has order $2j$, while
the sign-preserving cycle $(a_1, a_2, \cdots, a_j)$ has order $j$.

\begin{cor} 
Let $B_L^+ \subset B_L$ be the subgroup of $B_L$ that preserves the orientation of $\hat L_0$, where $(\hat L_0, \hat L_1, \cdots, \hat L_k)$ is a $(k+1)$-component hyperbolic link in $S^3$ such that $(\hat L_1, \cdots, \hat L_k)$ is the trivial link. Since $B_L^+$ acts on $\hat L_0$ by translations, it is a cyclic group. Let $n$ be the order of the cyclic group $B_L^+$. The representation $$B_L^+ \to \Sigma_k \wr \Zed_2$$
is conjugate to a product of (signed) cycles and there are at most $5$ different cycle types can realized in the cycle decomposition of this action. Using the conventions from Proposition \ref{maxsym}, the action of $B_L^+$ on $S^3$ is conjugate to a $(p,q)$-action for some pair of integers $(p,q)\in \Zed^2$ with $GCD(p,q)=1$.  The possible cycles that can appear can have:
\begin{enumerate}
\item length $n$, preserving sign. These correspond to components of $L$ which can be separated from the singular set of the action of $A$ on $S^3$ by round balls. 
\item If $GCD(q,n)>1$ cycles of length $n/GCD(q,n)$, preserving sign. These are represented by components of $L$ which Hopf
link the singular set $(\{0\} \times \mathbb C) \cap S^3$.
\item If $GCD(p,n)>1$ cycles of length $n/GCD(p,n)$, preserving sign.  These are represented by components of $L$ which Hopf link the singular set $(\mathbb C \times \{0\}) \cap S^3$. 
\item If $GCD(p,n)=2$, cycles of length $n/2$, reversing sign. These are represented by components of $L$ which intersect the singular set $(\mathbb C \times \{0\}) \cap S^3$ in two points. 
\item If $GCD(q,n)>1$, there can be a cycle of length $1$, preserving sign.  This corresponds to a single component of $\hat L_1 \cup \cdots \cup \hat L_k$ coinciding with a component of the singular set of the action, $(\{0\} \times \mathbb C) \cap S^3$. 
\end{enumerate}
Moreover, (5) and (2) are exclusive.  Thus if (5) holds, $k-1$ is a non-negative integer-linear combination of
$n$ and $n/GCD(p,n)$.  If (5) does not hold, $k$ is a non-negative integer-linear combination of
$n$, $n/GCD(q,n)$ and $n/GCD(p,n)$.
\end{cor}

\begin{eg} Sakuma's example where $B_L \simeq D_{10}$.  $B_L^+$ is cyclic of order $10$ acting on $S^3$ via a $(5,2)$-action. $B_L^+ \to \Sigma_5 \wr \Zed_2$ (with indicated orientations, taking the generator of $B_L^+$ to be counter-clockwise rotation in the plane of the figure by $2\pi/5$ and rotation by $\pi$ in the direction of the axis orthogonal to the plane) is the cycle $(1, 2, 3, 4, 5)-$. 
{
\psfrag{l0}[tl][tl][1][0]{$L_0$}
\psfrag{l1}[tl][tl][0.8][0]{$L_1$}
\psfrag{l2}[tl][tl][0.8][0]{$L_2$}
\psfrag{l3}[tl][tl][0.8][0]{$L_3$}
\psfrag{l4}[tl][tl][0.8][0]{$L_4$}
\psfrag{l5}[tl][tl][0.8][0]{$L_5$}
$$\includegraphics[width=6cm]{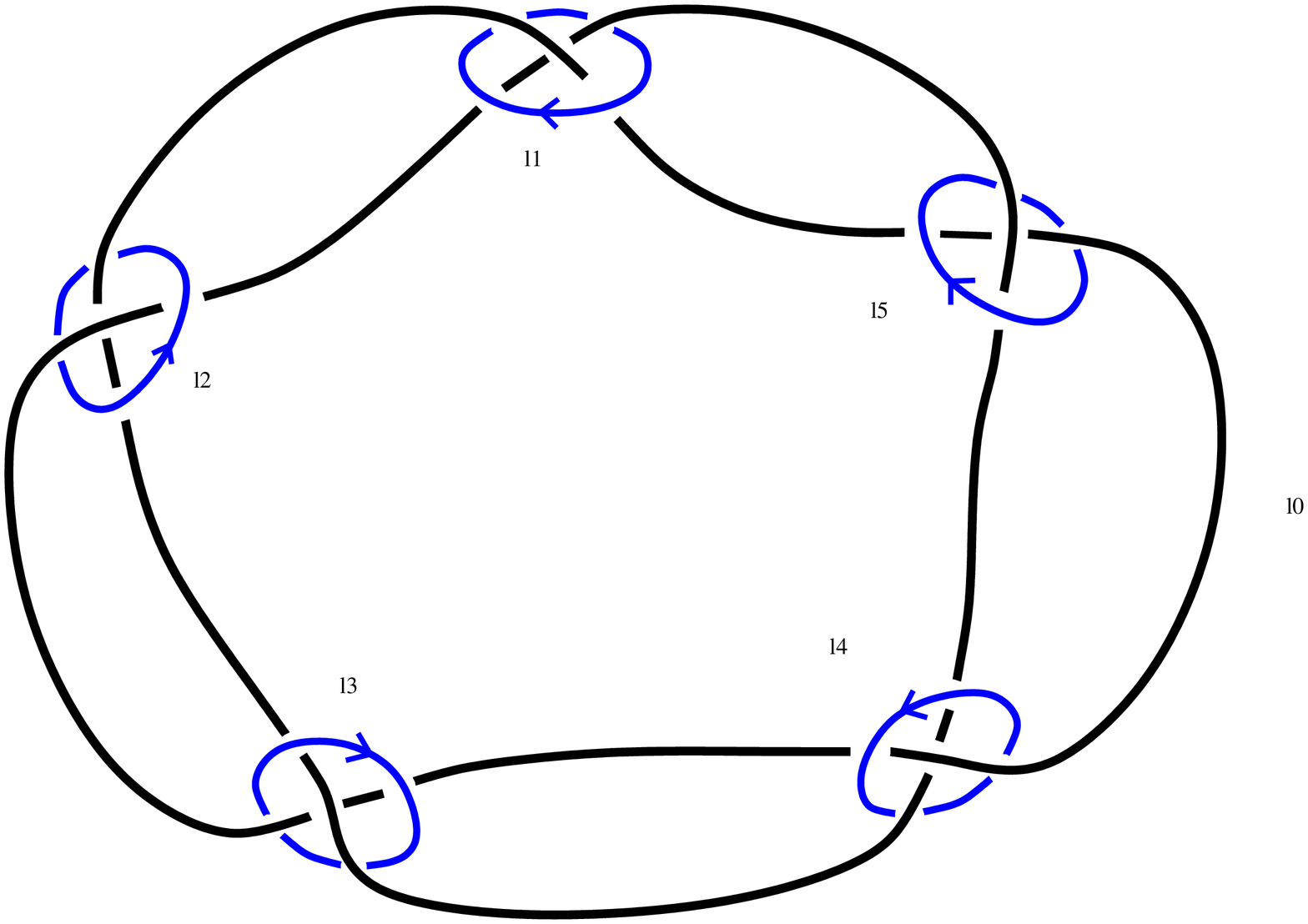}$$
}
\end{eg}

\begin{eg}
A hyperbolic example where $B_L = D_6$ giving a $(3,2)$-action on $S^3$.  Moreover one of the components of the link coincides with a singular circle of the action of $B_L$ on $S^3$. 
{
\psfrag{l0}[tl][tl][1][0]{$L_0$}
\psfrag{l1}[tl][tl][0.8][0]{$L_1$}
\psfrag{l2}[tl][tl][0.8][0]{$L_2$}
\psfrag{l3}[tl][tl][0.8][0]{$L_3$}
\psfrag{l4}[tl][tl][0.8][0]{$L_4$}
\psfrag{l5}[tl][tl][0.8][0]{$L_5$}
\psfrag{l6}[tl][tl][0.8][0]{$L_6$}
\psfrag{l7}[tl][tl][0.8][0]{$L_7$}
$$\includegraphics[width=6cm]{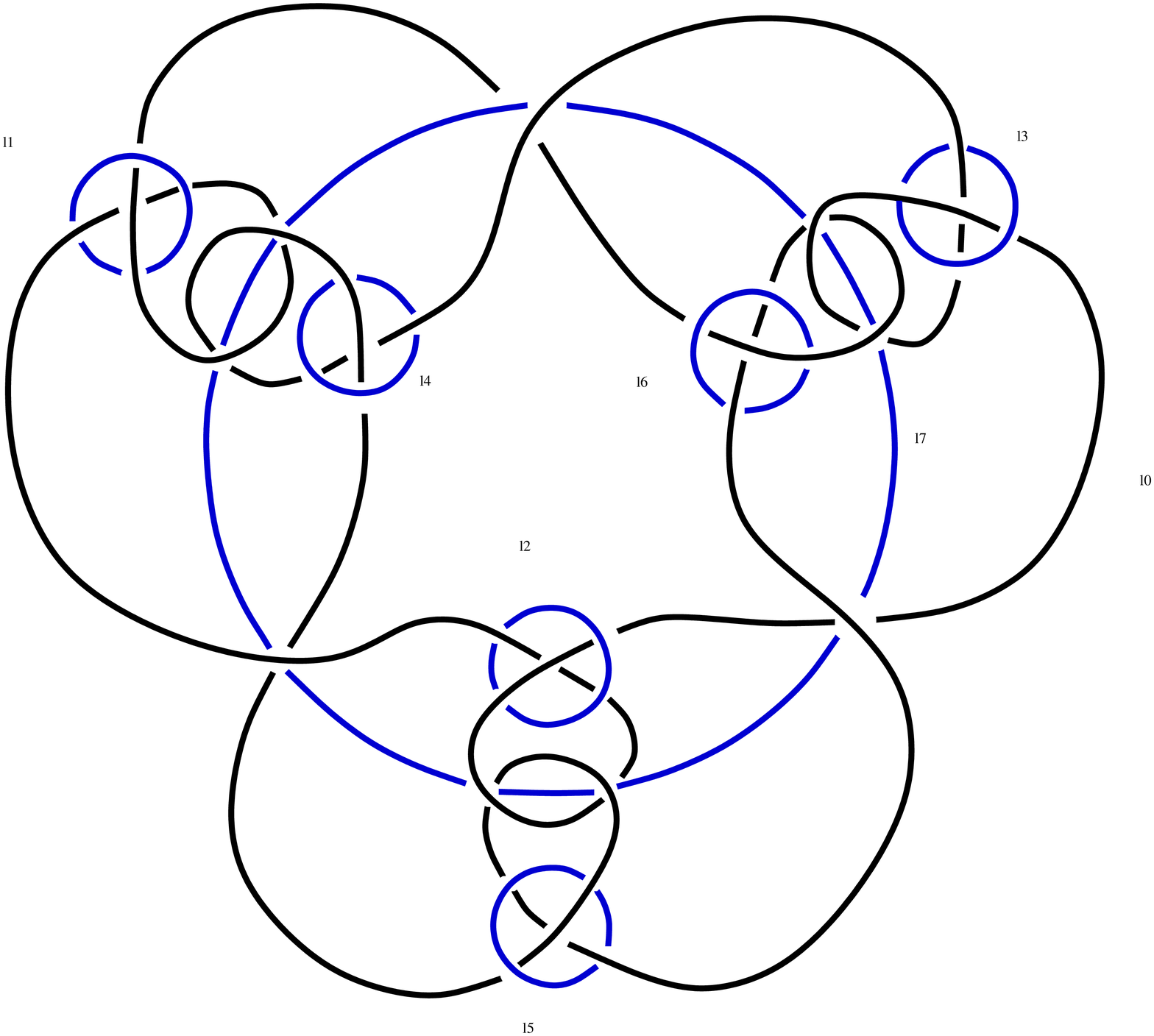}$$
}
$B_L^+$ is cyclic of order $6$, $B_L^+ \to \Sigma_7 \wr \Zed_2$ having cycle-type $(1,2,3,4,5,6)(7)$. 
\end{eg}

\section{Future directions}\label{futuredir}

This section points out some lines of inquiry that may be productive.  

\begin{prob}
Compute the homology of $\SP_{3,1}$ and $\SD_j^{D^{n-j}}$ {\it as an operad}.  
Does $\SD_1^{D^{n-1}}$ give any interesting homology operations on $H_* \EK{1,D^{n-1}}$ 
not provided by the $2$-cubes action on $\EK{1,D^{n-1}}$?
\end{prob}

For $H_* \SP_{3,1}$ a starting-point would be the work \cite{Cohen}. 

There is a wider class of embedding space that admits a `splicing operad' action.  Given a manifold $N$ with a co-dimension zero submanifold $V$, denote the space of embeddings $N \to N$ with support contained in $V$ by $Emb_V(N,N)$.  $\ED{j,M}$ would be the case $N = \Real^j \times M$ and $V = D^j \times M$.  Assume that $V$ is a manifold with co-dimension $2$ cubical corners.  Moreover, assume $\partial V$ is partitioned into two smooth manifolds with a common boundary $\partial V = W_1 \cup_C W_2$, $C$ the
co-dimension $2$ corner stratum.  We assume $W_1 \subset \partial N$ and $W_2$ is properly embedded in $N$. 
The associated operad to $Emb_V(N,N)$ would consist of equivalence classes $(k+2)$-tuples $(L_0, \cdots, L_k,\sigma)$ with $L_0 \in Emb_V(N,N)$ and $L_i: V \to V$ a self-embedding of $V$, just as in the definition of $\SD_j^M$.   Call this construction the operad of self-embeddings for the pair $(N,V)$.  Possibly interesting operads of this type would be when $N$ the total-space of a fibre bundle over a closed manifold ($p : N \to X$) with $V = p^{-1}(A)$, $A \subset X$ a co-dimension $0$ submanifold. 

\begin{prob}
Are operads of self-embeddings `interesting' outside of the $\SD_j^M$ and $\SC_j^M$ cases?  Do they fit into larger structures -- are there more general higher algebraic structures encoding the basic structure of diffeomorphism groups of manifolds? 
\end{prob}

The above problem is closely connected to a desire (shared by many) for spaces like $\K_{n,1}$ to have an action of the operad of framed $2$-discs, or some equivalent operad.  

An important difference between the descriptions of $\K_{3,1}$ as an algebra over the operads $\Cu_1'$ and $\SP_{3,1}$ respectively is that, although they are both free, the description of $\K_{3,1}$ over $\Cu_1'$ involves thinking of 
$\Cu_1'$ as an operad with non-empty `base' $\Cu_1'(0)$, while by design $\SP_{3,1}(0) = \emptyset$. The augmentation
maps for $\Cu_1'$ consist of `deleting an interval'.  $\SD_1^{D^2}(0)$ is non-empty, but in this case the augmentation
maps consist of operations including, among others, puck-deletion.  If one deletes a component of the Borromean rings, one gets a reducible link, thus $\SP_{3,1}$ has to have empty base.  To get puck-deletion as the augmentation map for $\SD_1^{D^2}$
one has to choose $Id_{\Real \times D^2}$ as the base-point of $\SD_1^{D^2}(0)$.  If one chooses a non-trivially framed
unknot, the augmentation maps become `twist maps' along the puck that is deleted.  Thus the operad $\SD_1^{D^2}$ encodes
a basic form of knot diagrammatics, and therefore must have a very rich homotopy-type, while $\SP_{3,1}$ is much more rigid
and has a relatively simple homotopy-type, as an operad. 

\begin{prob}
Is there a useful spaces-of-knots level description of further knot diagrammatics?  For example, is there an operadic or suitable higher-algebraic formalism for rational tangle decompositions of links \cite{BS}, or D.~Thurston's knotted trivalent graph constructions \cite{DT}? Further afield, perhaps the complexes describing spaces of connect-sum decompositions of manifolds \cite{CSR, HL} have an enlightening operadic formalism.   
\end{prob}

\providecommand{\bysame}{\leavevmode\hbox to3em{\hrulefill}\thinspace}


\end{document}